\documentclass[12pt, 14paper,reqno]{amsart}
\usepackage{amsmath,amsfonts,amssymb}
\usepackage{mathrsfs}
\usepackage{longtable}
\usepackage[breaklinks]{hyperref}
\usepackage{graphicx}
\usepackage{mathtools}
\usepackage{multirow}
\usepackage{multicol}
\usepackage{array}
\usepackage{tabularray}
\usepackage{enumitem}
\DeclarePairedDelimiter{\ceil}{\lceil}{\rceil}

\theoremstyle{plain}
\newcommand{\floor}[1]{\lfloor #1 \rfloor}
\newtheorem{thm}{Theorem}[section]
\newtheorem{lem}{Lemma}[section]

\newtheorem{cor}{Corollary}[section]

\newtheorem{prop}{Proposition}[section]

\newtheorem{thma}{Theorem}

\theoremstyle{proof}

\numberwithin{equation}{section}

\begin{document} 
\title[Sums of integral squares]{Sums of squares of  integer-multiple of an integral element on real bi-quadratic fields}
\author{Srijonee Shabnam Chaudhury}
\address{ Srijonee Shabnam Chaudhury
@Harish-Chandra Research Institute, A CI of Homi Bhabha National Institute, Chhatnag Road, Jhunsi,  Allahabad 211 019, India.}
\email{srijoneeshabnam@hri.res.in}
\maketitle
\begin{abstract}
 For any given positive integer $m$ we construct certain totally positive algebraic integers $\alpha$ of a real bi-quadratic field $K$ and obtain some necessary conditions for which $m\alpha$ can not be represented as sum of integral squares. We show this for integers lie in quadratic subfields of $K$ and for integers which are in $K$ but not in any quadratic subfield of $K$. We provide examples in tabular form for each cases to corroborate the results.
\end{abstract}
\maketitle
\section{Introduction}
The representation of an integer by sum of integer squares or by other quadratic forms in a number field $K$ has a long played fundamental role in algebraic number theory. It dates back at least to Pythagorean triplet, that is, the representation of $0$ by the ternary form $x^2 + y^2 - z^2$, from around $1800$ BC. The study of Pell's equation, that is, representation of $1$ and other small integers by the binary form $x^2 - dy^2$ (where $d$ is a square-free positive integer), was considered to be started in $400$ BC by Greek mathematicians and later was carried on by Archemides, Diophantus, Brahmagupta, Bhaskara and by many more authors.

The Modern history begins with stalwart mathematicians, like Fermat, Euler, Gauss who studied representation of primes by binary forms $x^2 + dy^2$ (where $d$ is a square-free positive integer). Finally, in 1770, Lagrange proved his famous four square theorem \cite{LA} which says that every positive integer can be represented by the positive definite quaternary form $x^2 + y^2 + z^2 + w^2$. At the beginning of twentieth century, Ramanujan \cite{R} and Dickson \cite{D} extended Lagrange's result and they showed that up-to equivalence, there exist only a finite number of positive definite quaternary forms which represent all positive integers. The most recent developments were made by Conway-Schneeberger \cite{Co}, Bhargaba- Hanke \cite{BH}, and Rouse \cite{Ro}, which gave equivalent conditions for a positive definite form to represent all positive integers and other special infinite sets.

Now, Let $K$ be a number field, $\mathcal{O}_K$ be its ring of integer and $\mathcal{O}_K^+(+)$ be the additive semigroup of totally positive integers. It is of considerable interest to look into the above mentioned results in $K$. In 1902, Hilbert conjectured that every totally positive element in $K$ can be written as a sum of four squares. The first published proof was by Siegel \cite{SI21} in 1921. He \cite{SI45} also proved in 1945 that $\mathbb{Q}$ and $\mathbb{Q}(\sqrt{5})$ are the only two totally real fields in which every element in $\mathcal{O}_K^+$ can be written as sum of integer squares. If $K$ is non totally real, Siegel proved that, every element in $\mathcal{O}_K^+$ can be written as sum of integer squares if and only if the discriminant of $K$ is odd. 

Now it is natural to ask whether for a given $m \in \mathbb{N}$ there are only finitely many totally real number fields for which $m$-th multiple of every integer can be written as sum of integral squares. For real quadratic field $\mathbb{Q} \sqrt{D}$ (where $D > 1$ is a square-free positive integer) Kala and Yatsyna \cite{KY} have given an upper bound of $m$ depending on $D$ such that if $m$ is less than that bound then all element of $m\mathcal{O}_K^+$ can not be written as sum of squares in $\mathcal{O}_K$. Later Raska \cite{Ra} improved this bound for quadratic number fields .
In this article, we find answer to similar questions in the set up of real bi-quadratic field $K = \mathbb{Q}(\sqrt{p}, \sqrt{q})$ where $p$ and $q$ are two distinct, positive and square-free integers $\geq 2$ and $p < q < r$ where $r = \frac{pq}{gcd(p,q)^2}$.


In Propositions \ref{prop1}, \ref{prop12} and \ref{prop61}, we construct certain intervals depending on positive integers $m$, $\kappa$ and $t$ and investigate the existence of $m$-th multiple of an algebraic integer in some quadratic subfield of $K$ which can not be represented as sum of squares in $\mathcal{O}_K$.

Then in Theorems \ref{thm1}\ref{thm1a}, \ref{thm2}\ref{thm2a} and \ref{thm3}\ref{thm3a} we group together each of these intervals around a fixed ratio of $\kappa : t$ and get large ones. This leads us to find a sharp lower bound of '$m$' in Corollaries \ref{cor2}\ref{cor2a},\ref{cor3}\ref{cor3a}, \ref{cor3b} and \ref{cor4}\ref{cor4a} respectively. 
As a Corollary of Theorem \ref{thm1}\ref{thm1a} we get:
\begin{thma}(Corollary \ref{cor2}\ref{cor2a})
    
Let $K = \mathbb{Q}(\sqrt{p},\sqrt{q}) $ where $p\equiv 2 \pmod 4$ , $q \equiv 3 \pmod 4$. Let  $m$, $t \in \mathbb{N}$ and $m \leq$ min $\{ \frac{r}{4 \ceil{\kappa \sqrt{r}}}, p, q \} $ for some positive integer $\kappa$ and $F_1 = \mathbb{Q}(\sqrt{p})$, $F_2 = \mathbb{Q}(\sqrt{q})$.  

If $p,q \geq \left(\frac{m}{2}+4 \right)^2$ then not all elements of $m\mathcal{O}_{F_1}^+$ and $m\mathcal{O}_{F_2}^+$ are represented as the sum of squares in $\mathcal{O}_K$.
\end{thma}
Corollaries \ref{cor3}\ref{cor3a}, \ref{cor3b} and \ref{cor4}\ref{cor4a} give similar lower bounds.

On the other hand, Propositions \ref{prop2}, \ref{prop22} and \ref{prop62} exhibit certain intervals depending on $m, \kappa, t_1, t_2 \in \mathbb{N}$ and some fixed relation between $p,q$ and $r$. Then uniting each of these intervals using the ratio of $\kappa:t_1$ and $\kappa:t_2$ in Theorem \ref{thm1}\ref{thm1b}, \ref{thm2}\ref{thm2b} and \ref{thm3}\ref{thm3b} we again get a lower bound of $m$ for which there exist integers in $m\mathcal{O}_K$ (but not in any quadratic subfield of $K$) which can not be represented as sum of integral squares.

Corollaries \ref{cor2}\ref{cor2b},\ref{cor3}\ref{cor3c} and \ref{cor4}\ref{cor4b} give lower bounds of these $m$'s for different types of integral elements in $K$ depending on $p,q \pmod 4$. As a sample, we state Corollary \ref{cor2}\ref{cor2b} where $p \equiv 2 \pmod 4$, $q \equiv 3 \pmod 4$.
\begin{thma}(Corollary \ref{cor2}\ref{cor2b})

Let $m < \frac{r}{16\kappa\ceil{\sqrt{r}}}$ for some positive integer $\kappa$. Let us also assume 
$$
N_1 = \text{max} \{\ceil{2\sqrt{q} - \sqrt{p}} , \ceil{\sqrt{r} - \sqrt{p}} \} -\sqrt{p} \}
$$
and
$$
N_2 = \text{max} \{  \ceil{2\sqrt{p} - \sqrt{q}} , \ceil{\sqrt{r} - \sqrt{q}} \} - \sqrt{q} \}
$$
are two constants.

If $p \geq \left( m + \sqrt{2mN_1}  \right)^2$ and $q \geq \left( m + \sqrt{2mN_2}  \right)^2$ then not all elements of $m\mathcal{O}_K^+$ are represented as the sum of squares in $\mathcal{O}_K$.
\end{thma}

Corollaries \ref{cor3}\ref{cor3c} and \ref{cor4}\ref{cor4b} give similar lower bounds.
\section{Some Basics}
\subsection{Totally positive integers and totally real field}
Let $K$ be an algebraic number field of degree $n$ with exactly $r_1$ real embeddings $ \sigma_1$,...,$\sigma_{r_1}$ and $r_2$ pairs of complex embeddings $\tau_1$,$\Bar{\tau_1}$,...,$\tau_{r_2}$,$\Bar{\tau_{r_2}}$, where $r_1+2r_2=n$. The field is totally real in the case $r_1=n$.\\
A number $\alpha$ in $K$ is called totally positive whenever the $r_1$ conjugates $\sigma_1(\alpha)$,..,$\sigma_{r_1}(\alpha)$ are all positive, denoted by $\alpha >> 0$. The set of all totally positive elements of $\mathcal{O}_K$ will be denoted by $\mathcal{O}_K^+$. This set is closed under addition and multiplication. We use the symbol $\alpha >> \beta$ to denote $\alpha - \beta >> 0$. Then the norm and trace of totally positive integers have the following basic properties:
\begin{itemize}
    \item  $Tr_{K/\mathbb{Q}}(\alpha) > 0 $.
    \item  $\mathcal{N}_{K/\mathbb{Q}}(\alpha) > 0$.
    \item  If $\alpha >> \beta$, then $Tr_{K/ \mathbb{Q}}(\alpha) >>  Tr_{K/ \mathbb{Q}}(\beta)$.
    \item  If $\alpha >> \beta$, then $\mathcal{N}_{K/ \mathbb{Q}}(\alpha) >>  \mathcal{N}_{K/ \mathbb{Q}}(\beta)$.
    \item  If $K$ is of degree $N$ over $\mathbb{Q}$ and $\alpha, \beta \in \mathcal{O}_K^+$, then by generalised Holder's inequality,
    $$
    \sqrt[N]{\mathcal{N}_{{K}/\mathbb{Q}}}(\alpha + \beta) >> \sqrt[N]{\mathcal{N}_{{K}/\mathbb{Q}}}(\alpha ) + \sqrt[N]{\mathcal{N}_{{K}/\mathbb{Q}}}(\beta).
    $$
Equality holds if and only if  $\frac{\alpha}{\beta}\in \mathbb{Q}$.

Here $\mathcal{N}_{K/\mathbb{Q}}\rightarrow$ is the norm operator on $K/ \mathbb{Q}$ and $Tr_{K/\mathbb{Q}}\rightarrow$ is the trace operator on $K/ \mathbb{Q}$.
\end{itemize}

\subsection{Bi-quadratic field and its integral basis}
Let $p$ and $q$ are two distinct, positive and square-free integers $\geq 2$, $g = gcd(p,q)$ and $r = \frac{pq}{g^2}$. Then $K= \mathbb{Q}(\sqrt{p},\sqrt{q})$ is called a real bi-quadratic field and $\mathbb{Q}(\sqrt{p})$, $\mathbb{Q}(\sqrt{q})$ and $\mathbb{Q}(\sqrt{r})$ are all three quadratic sub-fields of $K$ . Since $K$ is a splitting field of a polynomial $(x^2-p)(x^2-q)$ over $\mathbb{Q}$ therefore $K/\mathbb{Q}$ is a  degree $4$ Galois extension. Thus there are four embeddings of $K$ into $\mathbb{C}$: If $\alpha \in K$ is of the form $a+b\sqrt{p}+c\sqrt{q}+d\sqrt{r}$ with $a,b,c,d \in \mathbb{Q}$, the embeddings are:
\begin{eqnarray*}
\sigma_1(\alpha) =  a+b\sqrt{p}+c\sqrt{q}+d\sqrt{r}.\\
\sigma_2(\alpha) =  a-b\sqrt{p}+c\sqrt{q}-d\sqrt{r}.\\
\sigma_3(\alpha) =  a+b\sqrt{p}-c\sqrt{q}-d\sqrt{r}.\\
\sigma_4(\alpha) =a-b\sqrt{p}-c\sqrt{q}+d\sqrt{r}.
\end{eqnarray*}
Depending on $p$, $q \pmod4 $ and possibly interchanging the role of $p$,$q$ and $r$, every case can be converted into one of the following (see [\cite{Ja}, section 8]):
\begin{itemize}
\item[(1)] $ p \equiv 2 \pmod 4 \hspace{1cm} q\equiv 3 \pmod 4$.
\item[(2)] $ p \equiv 2 \pmod 4 \hspace{1cm}   q\equiv 1 \pmod 4$.
\item[(3)]$ p \equiv 3 \pmod 4 \hspace{1cm}   q\equiv 1 \pmod 4$.
\item [(4)] $ p \equiv 1 \pmod 4  \hspace{1cm}  q\equiv 1 \pmod 4$  and
\begin{itemize}
     \item [(i)] $\frac{p}{gcd(p,q)} \equiv \frac{q}{gcd(p,q)} \equiv 1 \pmod 4 $ or
     
         \item [(ii)] $ \frac{p}{gcd(p,q)} \equiv \frac{q}{gcd(p,q)} \equiv 3 \pmod4 $.
\end{itemize}

\end{itemize}
    
An integral basis of $\mathcal{O}_K$ in different cases has the following form (See [\cite{Wi}, Theorem 2]):
\begin{itemize}
    \item [(1)] $\mathcal{B}_1 = \{1, \sqrt{p}, \sqrt{q}, \frac{\sqrt{p}+\sqrt{r}}{2} \}$.
    \\
\item[(2)]  $\mathcal{B}_2 = \{1, \sqrt{p}, \frac{\sqrt{q}+1}{2}, \frac{\sqrt{p}+\sqrt{r}}{2} \}$.
\\
\item[(3)]  $\mathcal{B}_3 = \{1, \sqrt{p}, \frac{\sqrt{q}+1}{2}, \frac{\sqrt{p}+\sqrt{r}}{2} \}$. 
\\
\item [(4.a)] $\mathcal{B}_{41} = \{1, \frac{\sqrt{p}+1}{2}, \frac{\sqrt{q}+1}{2}, \frac{1+\sqrt{p}+\sqrt{q}+\sqrt{r}}{4} \}$.
\\
\item[(4.b)] $ \mathcal{B}_{42} = \{1, \frac{\sqrt{p}+1}{2}, \frac{\sqrt{q}+1}{2}, \frac{1-\sqrt{p}+\sqrt{q}+\sqrt{r}}{4} \}$.
\end{itemize}
In all these cases $p \equiv r \pmod4 $. Thus in the cases $1,2$ and $3$, $p$ and $r$ are interchangeable. In case $4$, all of $p$, $q$ and $r$ are interchangeable.
\subsection{Quadratic Forms:} 
Let $K$ be a number field of degree $N$ over $\mathbb{Q}$ and let $\{ x_1, x_2,..., x_d \}$ are $d$ variables. Then a $d$-ary quadratic form over $\mathcal{O}_{K}$ is the expression
$$
\mathbf{Q}(x_1,...,x_d) = \sum_{1\leq i\leq j\leq n} a_{ij}x_{i}x_{j}
$$
where $a_{ij} \in \mathcal{O}_K$.\\
Every quadratic form $\mathbf{Q}(x_1,...,x_d)$ can be associated with a symmetric matrix as follows: 
$$
\mathbf{Q}(x_1,..x_d) = \begin{pmatrix}x_1 & \cdots & x_d
\end{pmatrix} 
\begin{pmatrix}a_{1,1} & \frac{a_{1,2}}{2} & \cdots & \frac{a_{1,d}}{2} \\
\frac{a_{2,1}}{2} & a_{2,2} & \cdots & \frac{a_{2,d}}{2} \\
\vdots  & \vdots  & \ddots & \vdots  \\
\frac{a_{d,1}}{2} & \frac{a_{d,2}}{2} & \cdots & a_{d,d} 
\end{pmatrix} 
\begin{pmatrix} x_1 \\
\vdots &
\\
\vdots &
\\
x_d
\end{pmatrix}
$$
We denote the matrix by $\mathbf{A}$. This matrix is called the Gram matrix associated with $\mathbf{Q}$. Then $\mathbf{Q}$ is said to be
\begin{itemize}
    \item positive, definite if $\mathbf{A}$ is a positive, definite matrix,
    \item diagonal if $\mathbf{A}$ is a diagonal matrix, and,
    \item universal if it represents all totally positive integers of $K$.
\end{itemize}
Thus an element can be written as a sum of integral squares if and only if the Gram matrix associated with its quadratic form is the identity matrix. For more details on quadratic form see O'Mear's book \cite{OM}
\section{Case I}\label{c1}
In this section we discuss our results for the field $K = \mathbb{Q}(\sqrt{p},\sqrt{q}) $ where $p\equiv 2 \pmod 4$ and $q \equiv 3 \pmod 4$. Here the basis of $\mathcal{O}_K$ is --
$$
\mathcal{B}_1 = \{1, \sqrt{p}, \sqrt{q}, \frac{\sqrt{p}+\sqrt{r}}{2} \}.
$$
Let us consider $\alpha = u + \frac{\kappa_1}{2}\sqrt{p} + \kappa_2 \sqrt{q}+ \frac{\kappa_3}{2}\sqrt{r}\in \mathcal{O}_K^+$ where $\kappa_1$,$\kappa_2$,$\kappa_3$ are any three fixed non-negative integers. Now the question is, whether some integer multiple of $\alpha$ can be written as $\sum \alpha_i^2$, where $\alpha_i\in \mathcal{O}_K$ for all $i$. As an arbitrary element of $\mathcal{O}_K$, we can write $\alpha_i$'s as
$$
\alpha_i = \sum 
 a_i + \frac{b_i}{2}\sqrt{p}+ c_i\sqrt{q} + \frac{d_i}{2}\sqrt{r}
$$
where $a_i, b_i, c_i, d_i \in \mathbb{Z}$.

If $\alpha $ belongs to some quadratic subfield of $K$, that is, $\alpha$ can be written in the form $u + \kappa\sqrt{D}$ where $D \geq 2 $ is a square-free positive integer or half integer, then it is easy to find a minimal value of $u$ depending on $\kappa$ and $D$ so that $u + \kappa\sqrt{D}$ become the smallest totally positive integer of that form. But if not, then there exist no such $u$ depending on $\kappa_1$, $\kappa_2$ $\kappa_3$ for which $\alpha \in \mathcal{O}_K^+$ attains the minimal value. Still for finding the answer of this question it is important to find a suitable value of $u$ and to estimate a lower bound of the trace of $\sum \alpha_i^2$. If this lower bound is greater than $mu$ for some $m \in \mathbb{N}$ then we obtain $\alpha\in \mathcal{O}_K^+$ such that $m\alpha$ can not be written as sum of squares.

For this purpose we need two lemmas. One is for the quadratic case and other is for biquadratic fields. Although the proof of first one is similar to [\cite{Ra}, Lemma 3], we give the statement for the sake of completeness.
\begin{lem}\label{lem1}
Let $a_i$, $b_i$ be non- negative integers satisfying $\sum_{i} a_ib_i = \kappa $ for a fixed positive integer $\kappa$. Let $D$ be real number and $t$ be a positive integer. Then there exists an interval
 $$
I_t(\kappa)=\begin{cases} \left[ \frac{\kappa^2}{t(t+1)},\frac{\kappa^2}{t(t-1)} \right] & \text{if} \hspace{5mm} t >1, \\
\left[ \frac{\kappa^2}{2},\infty \right) & \text{if} \hspace{5mm} t = 1,
\end{cases}
$$
such that
\begin{itemize}
    \item [(i)] If $D \in I_t(\kappa) $ then $\sum_{i} a_i^2 + Db_i^2 \geq \frac{\kappa^2}{t}+ tD$.
    \item [(ii)] If $D \in I_t(2\kappa) $ then $\sum_{i} a_i^2 + D\frac{b_i^2}{4} \geq \frac{\kappa^2}{t}+ \frac{tD}{4}$.
\end{itemize}
\end{lem}
Let $f(x) = \frac{\kappa^2}{x} + xD$. Then from the above lemma we can say that minimum of $f(x)$ over all positive integers is equal to $f(t)$ if and only if $D \in I_t(\kappa)$. But for biquadratic fields this is not true.
Next lemma deals with the case for elements which belong to $K$ but not in any quadratic subfield of $K$.
\begin{lem}\label{lem2}
Let $a_i,b_i,c_i$,  $1\leq i \leq n$ be non-negative integers satisfying
$$
\sum a_ib_i = \kappa_1 \hspace{5mm} \sum a_ic_i = k_2 \hspace{5mm} \sum b_ic_i= \frac{\kappa_3}{g}
$$

where $\kappa_1, \kappa_2,\kappa_3$ are any three fixed positive integers.

If 
$$
p\in \left[ \frac{\kappa_1^2}{t_1^2}, \frac{\kappa_1^2}{t_1(t_1+1)}\right] = L_{t_1}(\kappa_1) \hspace{5mm} \text{and} \hspace{5mm} q\in \left[ \frac{\kappa_2^2}{t_2^2}, \frac{\kappa_2^2}{t_2(t_2+1)}\right] = L'_{t_2}(\kappa_2)  
$$
where $t_1, t_2$ are any two positive integers.

Then,
$$
\sum a_i^2+ b_i^2p + c_i^2q \geq \frac{\kappa^2}{t}+ t_1p + t_2q
$$
where $t=$ min$\{ t_1,t_2 \}$ and $\kappa=$ max $\{\kappa_1, \kappa_2 \}$.
\end{lem}
\begin{proof}
Let $l_1= \sum b_i^2$ and $l_2= \sum c_i^2$. Then by Cauchy-Schwarz inequality
$$
\sum a_i^2 \geq \frac{\sum (a_ib_i)^2}{\sum b_i^2}= \frac{\kappa_1^2}{l_1} \hspace{5mm} \text{and} \hspace{5mm} \sum a_i^2 \geq \frac{\sum (a_ic_i)^2}{\sum c_i^2}= \frac{\kappa_2^2}{l_2}.
$$
This implies,
$$
\sum a_i^2+ b_i^2p + c_i^2q \geq \frac{\kappa^2}{l}+ l_1p + l_2q,
$$
where $l=$ min $\{l_1, l_2\}$ and $\kappa = \text{max} \{\kappa_1, \kappa_2 \}$.

Therefore, we have to show that,
\begin{equation}\label{eqn1}
\frac{\kappa_1^2}{l_1}+ l_1p+ l_2q \geq \frac{\kappa_1^2}{t_1}+ t_1p + t_2q,
\end{equation}
and
\begin{equation}\label{eqn2}
\frac{\kappa_2^2}{l_2}+ l_1p+ l_2q \geq \frac{\kappa_2^2}{t_2}+ t_1 p+ t_2q.
\end{equation}
\eqref{eqn1} implies $(l_1-t_1)(l_1t_1p - \kappa_1^2) + q(l_2-t_2) \geq 0$,

which is true if 
\begin{equation}\label{eq3}
t_1(t_1-1)p \leq \kappa_1^2 \leq t_1(t_1+1) \hspace{5mm} \text{and} \hspace{5mm} l_2 \geq t_2.
\end{equation}
Similarly, \eqref{eqn2} implies $(l_2-t_2)(l_2t_2q - \kappa_2^2) + p(l_1-t_1)$,

which is true if 
\begin{equation}\label{eq4}
   t_2(t_2-1)q \leq \kappa_2^2 \leq t_2(t_2+1) \hspace{5mm} \text{and} \hspace{5mm} l_1 \geq t_1.
\end{equation}
Both \eqref{eq3} and \eqref{eq4} are true since $p \in L_{t_1}(\kappa_1) $ and $q \in L'_{t_2}(\kappa_2) $.
\end{proof}
In the above lemma, the minimum of the function $g(x,y) = \frac{\kappa^2}{x} + x p + yq$ over positive integers is $g(t_1,t_2)$ if $p\in  I_{t_1}(\kappa_1)$  and $q\in I'_{t_2}(\kappa_2)$. But this is not the 'only if' condition for attaining this minimum value.

On the other hand, in this lemma we get equality if 
$$
a_i= r_1b_i \hspace{5mm} a_i=r_2c_i \hspace{5mm} \text{for some fixed real numbers $r_1$ and $r_2$}
$$
$$
l_1= \sum b_i^2 = t_1 \hspace{5mm} \text{and} \hspace{5mm} l_2= \sum c_i^2 = t_2.
$$
From these conditions we get 
$$
r_1t_1 = \kappa_1 \hspace{5mm} \text{and} \hspace{5mm} r_2t_2 = \kappa_2.
$$
These implies that,
$$
\frac{b_i \kappa_1}{t_1} = \frac{c_i \kappa_2}{t_2} = a_i.
$$
If $d_1 =$ gcd $(\kappa_1, t_1)$ and $d_2 =$ gcd $(\kappa_2, t_2)$ then $t_1| b_i\kappa_1 \implies t_1' | b_i \kappa_1'$ where $t_1' = \frac{t_1}{d_1}$. This implies $t_1'^2| t_1$. Similarly, $t_2'^2| t_2$ where $t_2' = \frac{t_2}{d_2}$. Conversely, if $t_1'^2| t_1$, $t_2'^2| t_2$ and $\frac{\kappa_1}{\kappa_2} = \frac{d_1}{d_2}$ then we find suitable $a_i,b_i,c_i$ for which the equality holds. 

For example, if we take $b_i = t_1'$, $c_i = t_2'$ then, $a_i = r_1b_i = \frac{\kappa_1}{t_1}t_1'= \kappa_1'$ and similarly $a_i = \kappa_2'$, where $\kappa_1' = \frac{\kappa_1}{d_1}$ and  $\kappa_2' = \frac{\kappa_2}{d_2}$. Since they are equal so it is possible and the equality holds for all $p,q$.

As mentioned above if $\alpha = u + \frac{\kappa}{2}\sqrt{D}$ lies in some quadratic subfield of $K$ then taking $u = \ceil{\frac{\kappa}{2}\sqrt{D}}$ we get the smallest totally positive integer of this form and in this case $\alpha' < 1 $.

Further we consider some integer multiples of $\alpha$ and using the above lemma will find when it can not be written as sum of squares of elements of $\mathcal{O}_K$.

\begin{prop}\label{prop1}
Let $m,t\in \mathbb{N}$. Let
$$
I_1= \left[\frac{m\kappa}{2t} + \frac{{\sqrt{m}}}{\sqrt{t}},\frac{m \kappa}{2(t-1)} - \frac{{\sqrt{m}}}{\sqrt{t-1}} \right]
$$
for $t > 1$ be an interval.

If $\kappa \in \mathbb{Z}_{>0}, m \leq$ min $\{ \frac{r}{4 \ceil{\kappa \sqrt{r}}}, p, q \} $ with $\sqrt{p}, \sqrt{q} \in I_{1}$, then not all the elements of $m\mathcal{O}_K^{+}$ are represented as sum of integral squares. These intervals are non-empty for 
    $$
    m\geq \frac{4t(t-1)(2t-1+2\sqrt{t(t-1))}}{k^2}.
    $$
    For $t=1$ we define the right bound of the intervals to be $\infty$ and are therefore always non-empty.

\end{prop}
\begin{proof}
Any $\alpha \in \mathcal{O}_K$ 
can be written as 
$$
\alpha = a + \frac{b}{2}\sqrt{p} + c\sqrt{q} + \frac{d}{2}\sqrt{r}
$$
 with $a,b,c,d \in \mathbb{Z}$  and $ 2|b-d $.
 
For arbitrary positive integer $\kappa$ consider $ \alpha = \ceil{\kappa\sqrt{p}} +\kappa\sqrt{p}$ and suppose $m \alpha = \sum_i \alpha_i^2 $ for some $\alpha_i = a _i+ \frac{b_i}{2}\sqrt{p} + c_i\sqrt{q} + \frac{d_i}{2}\sqrt{r} $ where $b_i \equiv d_i \pmod 2$. Now, comparing trace of $m\alpha$ with that of $\sum_i \alpha_i^2 $ we get 

$$
m \ceil{\kappa\sqrt{p}} = a _i^2 + \frac{b_i^2}{4} p  + c_i^2 q + \frac{d_i^2}{4} r.
$$
This is so as $ m < \frac{r}{4 \ceil{\kappa \sqrt{r}}} $ therefore $m \ceil{\kappa\sqrt{p}} <  \frac{r}{4}$.

Thus $d_i = 0 $ for all $i$ and then $2|b_i$. Therefore we substitute $b_i$ by $2b_i$ and get $\alpha_i = a _i+ b_i\sqrt{p} + c_i\sqrt{q} $ for all $i$.

Without loss of generality, we can choose $a_i \geq 0$ for all $i$ and $b_i \geq 0$ if $a_i=0$. Otherwise suppose there exist some $i$ such that $b_i < 0 , c_i < 0 $ and $a_i > 0$.
Then
$$
m \sigma_4(\alpha) \geq \sigma_4(\alpha_i)^2 = a_i^2 + b_i^2 p+ c_i^2q + 2a_i(-b_i) \sqrt{p} + 2a_i(-c_i) \sqrt{q} + 2b_ic_ig \sqrt{r} \geq 1 + p + 2\sqrt{p} \geq m,
$$
for $p \geq m $, where $g= gcd (p,q)$. Which is impossible because $\sigma_4 (\alpha) < 1$. Therefore both  $b_i$ and $c_i$ can not be negative for all $i$.

Further if possible let there exist $\alpha_i$ with $b_i < 0$ and $c_i > 0$. Then using $\sigma_2 (\alpha)$ instead of $\sigma_4 (\alpha)$ we get similar type of contradiction. Therefore we can assume $b_i \geq 0$ for all $i$.

Comparing the irrational parts in the expression of $m \alpha$, we get $m \kappa = 2\sum_i a_i b_i$ and that of rational parts from Lemma \ref{lem1} gives,
$$ 
m \kappa\sqrt{p} + m > m (\kappa\floor{\sqrt{p}} + 1 ) = \sum_i a_i^2 + b_i^2 p + c_i^2 q \geq a_i^2 + b_i^2 p \geq \frac{m^2\kappa^2}{t}+ tp
$$
for $p \in I_{t}(\frac{m\kappa}{2})$.
This is impossible if $\sqrt{p} \geq \frac{m \kappa}{2t}+ \frac{\sqrt{m}}{\sqrt{t}}$ or $\sqrt{p} \leq \frac{m \kappa}{2t} - \frac{\sqrt{m}}{\sqrt{t}}$.

Combining these constraints for $p \in I_{t}(\frac{m \kappa}{2})$ and $p \in I_{t-1}(\frac{m \kappa}{2})$ entails: if 
$$
\sqrt{p} \in \left[ \frac{m \kappa}{2t} + \frac{{\sqrt{m}}}{\sqrt{t}},\frac{m \kappa}{2(t-1)} - \frac{{\sqrt{m}}}{\sqrt{t-1}} \right] = I_1
$$
then $m\alpha$ can not be represented as a sum of squares.

This interval is non-empty only when  $|I_1| > 0$ and this implies 
$$
m \geq \frac{4t(t-1)(2t-1+ 2 \sqrt{t(t-1)}}{\kappa^2}
$$
Proceeding analogously we get $\sqrt{q} \in I_1$.
$I_1$ is well-defined for $t > 1$. For $t = 1$ we get the following interval
$$
\sqrt{p}, \sqrt{q} \in \left[\frac{m\kappa}{2} + \sqrt{m}, \infty \right) \hspace{5mm} \text{for} \hspace{5mm} p,q \in I_1(\frac{m\kappa}{2}).
$$
\end{proof}
An easy consequence of the above proposition shows that if two elements of the given form from two distinct subfields of $K$ can not be written as a sum of integral squares then their sum can not be represented as a square in $\mathcal{O}_K$. In some special cases it can not also be written as a sum of two squares in $\mathcal{O}_K$.
\begin{cor}\label{cor1}
    Let $\alpha_1 = \ceil{\kappa\sqrt{p}} +\kappa\sqrt{p} $ , $\alpha_2 = \ceil{\kappa\sqrt{q}} +\kappa\sqrt{q} $ and $\alpha_3 = \ceil{\kappa\sqrt{r}} +\kappa\sqrt{r} $   where $p, q,r$ are as defined above and $m$ is defined in proposition \ref{prop1}. If $m$ is even then 
$$
\frac{m}{2}(\alpha_i + \alpha_j)
$$
where $i,j \in \{1,2,3 \}$ and $i \neq j$,
can not be represented as a square in $\mathcal{O}_K$.

If $m =2$ and $\kappa = 1$ then also it can not be written as sum of two integral squares.
\end{cor}
\begin{proof}
    Let 
   $$
\frac{m}{2}(\alpha_1 + \alpha_2) = (a + \frac{b}{2}\sqrt{p} + c\sqrt{q} + \frac{d}{2} \sqrt{r})^2
$$
where $b \equiv d \pmod 2$.   

Proposition \ref{prop1} entails that $m \leq$ min $\{ \frac{r}{4 \ceil{\kappa \sqrt{q}}}, p, q \} $ and $p < q < r$. This implies $\frac{m}{2} < \frac{r}{4 (\ceil{\kappa \sqrt{p}} + \ceil{\kappa \sqrt{q}})} < \frac{r}{4} \implies \frac{m}{2}(\text{ rational part of $(\alpha_1 + \alpha_2)$) } < \frac{r}{4} \implies d = 0$. Therefore $b$ is even. Substituting $b$ by $2b$ gives
\begin{equation}\label{eqn111}
\frac{m}{2}(\alpha_1 + \alpha_2) = (a + b\sqrt{p} + c\sqrt{q} )^2.
\end{equation}
Comparing the irrational part of \eqref{eqn111} entails
\begin{equation}\label{eqn222}
    4ab = m\kappa
\end{equation}
\begin{equation}\label{eqn333}
    4ac = m\kappa
\end{equation}
and 
\begin{equation}\label{eqn444}
    4gbc = 0
\end{equation}
where $g = gcd(p,q)$.
As $g \neq 0$ from  \eqref{eqn444} entails either $b=0$ or $c=0$. In both the scenarios from  \eqref{eqn222} and \eqref{eqn333} we get that $m\kappa =0 $, which is impossible. Hence the first part is proved.

Let $m =2, \kappa = 1$ and 
\begin{equation}\label{eqn1111}
\frac{m}{2}(\alpha_1 + \alpha_2) = \frac{m}{2}\alpha = (\ceil{\sqrt{p}} + \ceil{\sqrt{q}} + \sqrt{p} + \sqrt{q}) = \sum_{i=1,2} (a_i + b_i\sqrt{p} + c_i\sqrt{q} )^2.
\end{equation}
Comparing irrational part in \eqref{eqn1111} entail
\begin{equation}\label{eqn2222}
a_1b_1 + a_2b_2 = 1
\end{equation}
\begin{equation}\label{eqn3333}
    a_1c_1 + a_2c_2 = 1
\end{equation}
and
\begin{equation}\label{eqn4444}
    c_1b_1 + c_2b_2 = 0.
\end{equation}
Further without loss of generality we may assume that $a_i \geq 0$. Let $b_i < 0$ and $c_i < 0$ for some $i$.Then
\begin{eqnarray*}
m \sigma_4(\alpha) &=& (a_i - b_i\sqrt{p} - c_i\sqrt{q})^2\\ 
&=& a_i^2 + b_i^2p + c_i^2q + 2a_ib_i\sqrt{p} + 2a_ic_i\sqrt{q} + 2gb_ic_i \sqrt{r}\\
&>& (p+q) > 2p > 2m
\end{eqnarray*}
implies that $\sigma_4 (\alpha) > 2$, which is not possible because 
$$ 
\sigma_4 (\alpha) = (\ceil{\sqrt{p}} + \ceil{\sqrt{q}} - \sqrt{p} - \sqrt{q}) < 2.
$$
Therefore both $b_i$ and $c_i$ can not be negative for any $i$. 
Thus we have two cases to be dealt with.\\
\textbf{Case I}:
$b_1c_1 \neq 0 \implies b_2c_2 \neq 0$ (from \eqref{eqn4444}).

If $b_1c_1 > 0$ and $b_2c_2 < 0$ then $b_1 > 0, c_1 > 0$ and either $b_2 < 0$ or $c_2 < 0$.

First we consider  $b_1, c_1, c_2 > 0$ and $b_2 < 0$. Since $a_1, a_2 \geq 0$ then $a_1c_1 \geq 0$ and $a_2c_2 \geq 0$. But from \eqref{eqn3333} we get both of them can not be positive. Therefore first consider $a_1 = 0$. This implies $a_2 = c_2 = 1$ and $a_2 = b_2 = 1$ (from \eqref{eqn2222}). But by assumption $b_2 < 0$.

Now let $a_2 = 0$. Then from \eqref{eqn2222} and \eqref{eqn3333}, we get $a_1 = b_1 = c_1 = 1$. Also \eqref{eqn4444} gives $b_2 = -1$ and $c_2 = 1$.

Putting these values in \eqref{eqn1111} we get 
$$
(\ceil{\sqrt{p}} + \ceil{\sqrt{q}} + \sqrt{p} + \sqrt{q}) = (1 + \sqrt{p} + \sqrt{q})^2 + (-\sqrt{p} + \sqrt{q})^2 = 1 + 2p + 2q + 2\sqrt{p} + 2\sqrt{q}
$$
which is not possible.

Now considering $b_1, c_1, b_2 > 0$ and $c_2 < 0$,  we get $a_1b_1 \geq 0$ and $a_2b_2 \geq 0$ which implies that either $a_1 = 0$ or $a_2 = 0$. Taking $a_1 = 0$ we get first the type of contradiction while in case when $a_2 = 0$, we get the second type of contradiction as above.\\
\textbf{Case II}:
$b_1c_1 = b_2c_2 = 0$

If $c_1 = c_2 = 0$, contradiction comes from $\eqref{eqn3333}$ and 
if $b_1 = b_2 = 0$, contradiction is gotten from $\eqref{eqn2222}$.

If $b_1 = 0 = c_2$ from \eqref{eqn3333} we get $a_2 = b_2 = 1$ and from equation \eqref{eqn2222} we get $a_1 = c_1 = 1$.

This implies

$$
(\ceil{\sqrt{p}} + \ceil{\sqrt{q}} + \sqrt{p} + \sqrt{q}) = (1 + \sqrt{p})^2 + (1 + \sqrt{q})^2 \implies p+q < 0
$$
which is not possible.
When $b_1 = 0 = c_2$, similar arguments will give the contradiction.

\end{proof}

If $\alpha = u + \frac{\kappa_1}{2}\sqrt{p} + \kappa_2 \sqrt{q}+ \frac{\kappa_3}{2}\sqrt{r}$ then it is not possible to find some uniform value of $u$ for which $\alpha$ is the smallest totally positive integer of this form. Also, we can not find any fixed real number, say $n$ such that conjugates of $\alpha$ is always less than that. Therefore for dealing with this case we have to consider some relation between $p, q$ and $r$.
\begin{prop}\label{prop2}
Let  $ \sqrt p + \sqrt{q} + 1 > \ceil{\sqrt{r}}$ and $m,t_1,t_2 \in \mathbb{N}$. Let $ \kappa$ be any arbitrary but fixed positive integer and $m < \frac{r}{16\kappa\ceil{\sqrt{r}}}$. Let
$$
C_1 = \text{max} \{\ceil{2\sqrt{q} - \sqrt{p}} , \ceil{\sqrt{r} - \sqrt{p}} \} -\sqrt{p} \},
$$
and
$$
C_2 = \text{max} \{  \ceil{2\sqrt{p} - \sqrt{q}} , \ceil{\sqrt{r} - \sqrt{q}} \} - \sqrt{q} \}
$$
be two constants.
Let 
$$
I_1' = \left[ \frac{m \kappa}{t_1} + \frac{{\sqrt{2m\kappa N_1}}}{\sqrt{t_1}},\frac{m \kappa}{(t-1)} - \frac{{\sqrt{2m\kappa N_1}}}{\sqrt{t_1-1}} \right] 
$$
for $t_1 >1$ 

$$
I'_2= \left[\frac{m\kappa}{t_2} + \frac{{\sqrt{2m\kappa N_2}}}{\sqrt{t_2}},\frac{m \kappa}{(t_2-1)} - \frac{{\sqrt{2m \kappa N_2}}}{\sqrt{t_2-1}} \right]
$$
for $t_2 > 1$ be two intervals.
Here $N_1$ and $N_2$ are any two constants depending on $C_1$ and $C_2$ respectively.

If $\sqrt{p} \in I'_{1}$ and $\sqrt{q}\in I'_{2}$ then not all elements of $m\mathcal{O}_K^{+}$ are represented as sum of integral squares. These intervals are non-empty for 
    $$
    m\geq \text{max} \{ \frac{2N_it_i(t_i-1)(2t_i-1+2\sqrt{t_i(t_i-1))}}{k} | i= 1,2\}.
    $$
    For $t_1 = t_2 = 1$ we define the right bound of the intervals to be $\infty$ and are therefore always non-empty.
\end{prop}  
\begin{proof}
    Let for an arbitrary positive integer  $\kappa$
    $$
    \alpha = u + 2\kappa \sqrt{p} + 2\kappa \sqrt{q} + \kappa \sqrt{r}
    $$ 
    where 
    $$
    u = 2 ~\mbox{max}~\{ \ceil{2\kappa \sqrt{p}}, \ceil{2\kappa \sqrt{q}}, \ceil{\kappa \sqrt{r}}\}.
    $$
Also  suppose 
$$
m \alpha = \sum_i \alpha_i^2  ~\mbox{for some}~ \alpha_i = a _i+ \frac{b_i}{2}\sqrt{p} + c_i\sqrt{q} + \frac{d_i}{2}\sqrt{r}, 
$$ 
where $b_i \equiv d_i \pmod 2$. As $m < \frac{r}{16 \kappa \ceil{\sqrt{r}}} \implies mu < \frac{r}{4} \implies d_i = 0$ for all $i$. This also implies, $b_i \equiv 0 \pmod 2$. Therefore we substitute $b_i$ by $2b_i$ and get $\alpha_i = a _i+ b_i\sqrt{p} + c_i\sqrt{q} $ for all $i$.
    
Without loss of generality, we can choose $a_i \geq 0$ for all $i$ and $b_i \geq 0$ if $a_i=0$. Otherwise, suppose there exist some $i$ such that $b_i < 0 , c_i < 0 $ and $a_i > 0$. Then 
\begin{eqnarray*}
m \sigma_4(\alpha) &\geq& 
 \sigma_4(\alpha_i)^2 \\
&=& a_i^2 + b_i^2 p + 2a_i(-b_i) \sqrt{p} + 2a_i(-c_i) \sqrt{q} + 2b_ic_ig \sqrt{r} 
\\
&\geq& 1+p + q + 2\sqrt{p} + 2\sqrt{q} + 2g\sqrt{r} \geq \ceil{\sqrt{r}}^2 \\
&\geq& 16m\kappa\ceil{\sqrt{r}}.
\end{eqnarray*}
This is an impossibility as $\sigma_4(\alpha) < 5\kappa \ceil{\sqrt{r}}$. 

Now, let us suppose that there exist $\alpha_i$ with $b_i < 0$ and $c_i > 0$. Then using $\sigma_2 (\alpha)$ instead of $\sigma_4 (\alpha)$ analogously we get contradiction. So $b_i \geq 0$ for all $i$.

The case of $\alpha_i $ with $b_i > 0$ and $c_i < 0 $ we can use $\sigma_3 (\alpha)$ (instead of $\sigma_4 (\alpha)$) and get similar contradiction. Thus $c_i \geq 0$ for all $i$. 

Comparing the irrational parts in the expression of $m \alpha$, we get 
$$
m \kappa = \sum_i a_i b_i = \sum_i a_i c_i = 2g\sum_i b_ic_i,
$$
and comparing rational parts from Lemma \ref{lem2} entails

$$
mu = \sum_i a_i^2 + b_i^2 p + c_i^2 q \geq \frac{m^2\kappa^2}{t_1}+ t_1p + t_2q.
$$
Taking $N_1 = C_1$ and comparing rational parts from Lemma \ref{lem1} we get 
\begin{eqnarray*}
2m \kappa(\sqrt{p} + N_1) > mu &=& \sum_i a_i^2 + b_i^2 p + c_i^2 q\\
&\geq& \frac{m^2\kappa^2}{t_1}+ t_1p + t_2q\\ 
&\geq& \frac{m^2\kappa^2}{t_1}+ t_1p,
\end{eqnarray*}
for $p \in L_{t_1}\left( m\kappa \right)$.
Combining these restrictions for $p \in L_{t_1}(m \kappa)$ and $p \in L_{t_1-1}(m \kappa)$ we derive that, if 
$$
\sqrt{p} \in \left[ \frac{m \kappa}{t_1} + \frac{{\sqrt{2m \kappa N_1}}}{\sqrt{t_1}},\frac{m \kappa}{(t_1-1)} - \frac{{\sqrt{2m \kappa N_1}}}{\sqrt{t_1-1}} \right] = I'_1
$$
then $m\alpha$ can not be represented as a sum of squares.
 
Proceeding similarly and replacing $ N_1, t_1$ by $ N_2, t_2$ simultaneously where $N_2 = C_2$, we get that $\sqrt{q} \in I_2'$. Further $I_1'$ and $I_2'$ are non-empty if 
$$
m\geq \text{max} \{ \frac{2N_it_i(t_i-1)(2t_i-1+2\sqrt{t_i(t_i-1))}}{\kappa} | i= 1,2\}.
$$
Also $I_1'$ and $I_2'$ are well-defined for $t_1, t_2 > 1$ and for $t_1 = t_2 = 1$  we get the following intervals:
$$
\sqrt{p} \in \left[m\kappa + \sqrt{2m\kappa N_1}, \infty \right) \hspace{5mm} \text{for} \hspace{5mm} p \in L_1(m\kappa)
$$
and
$$
\sqrt{q} \in \left[m\kappa+ \sqrt{2m\kappa N_2}, \infty \right) \hspace{5mm} \text{for} \hspace{5mm} p \in L_1'(m\kappa).
$$
\end{proof}
The results in Proposition \ref{prop1} and \ref{prop2} basically obtain some intervals depending on $m, \kappa$ and $t$ or $m, \kappa, t_i$ for $i= 1,2$. Now, using a fixed value of $\frac{\kappa}{t}$ or $\frac{\kappa}{t_i}$ for $i= 1,2$ respectively it is possible to construct large intervals depending on $m$, for which the previous results are true. We show that for this purpose $t$ must be finite and its upper bound depends on $m$. Therefore for a fixed value of $\kappa$ we always get finitely many intervals to unite.
\begin{thm}\label{thm1}
Let $m$ be a positive integer and $F_1 = \mathbb{Q}(\sqrt{p})$, $F_2 = \mathbb{Q}(\sqrt{q})$ where $p,q,m$ are as in Propositions \ref{prop1} and \ref{prop2}.
\begin{enumerate}
\item If $\sqrt{p}$, $\sqrt{q} \in \left[ \frac{m}{2}+4, \infty \right)$, $\left[ \frac{m}{2i}+ i\sqrt{40}, \frac{m}{2(i-1)} - (i-1)\sqrt{70} \right]$ for integer $i > 1$, then there exist elements in $m\mathcal{O}_{F_1}^{+}$ and $m\mathcal{O}_{F_2}^{+}$ which can not be written as sum of squares in $\mathcal{O}_K$. \label{thm1a}

\item If $\sqrt{p} \in \left[ m + \sqrt{2m N_1}, \infty \right)$, $ \left[ \frac{m}{i} + \sqrt{\frac{2mN_1}{i}}, \frac{m}{i-1}- \sqrt{\frac{2mN_1}{(i-1)}}\right]$ for integer $i > 1$ and $\sqrt{q} \in \left[ m + \sqrt{2m N_2}, \infty \right)$, $ \left[ \frac{m}{i} + \sqrt{\frac{2mN_2}{i}}, \frac{m}{i-1}- \sqrt{\frac{2mN_2}{(i-1)}}\right]$ for integer $i > 1$,\label{thm1b} 

then there exist elements of $m\mathcal{O}_K^+$ (but not elements of any quadratic subfield of it) which can not be written as sum of squares in $\mathcal{O}_K$.
\end{enumerate}
\end{thm}
\begin{proof}

 \begin{enumerate}

     \item See [\cite{Ra},Theorem $5$].

     \item  The technique is similar to part \eqref{thm1a}, but as the ratio of $(\kappa : t) $ appear in both the summands of the lower and upper bounds in given intervals, the calculation is far easier here. We include a proof of the result for $p$ for the sake of completeness.

Let $S(t , \kappa) = \left[ \frac{m \kappa}{t} + \frac{{\sqrt{2m \kappa N_1}}}{\sqrt{t}},\frac{m \kappa}{(t-1)} - \frac{{\sqrt{2m \kappa N_1}}}{\sqrt{t-1}} \right]$.

For the purpose of taking union of such intervals for fixed $m$, we make use of the ratio of $t$ and $\kappa$ in $\left[ \frac{m}{i}, \frac{m}{i-1} \right] = U(m,i)$ containing $\sqrt{p}$.

Observe that, when $\frac{t}{\kappa} \leq 1$ then $S(t, \kappa) \subseteq U(m,1) = \left[m, \infty \right)$. But for $\frac{t}{\kappa} < 1$, $S(t_1,\kappa_1) \subseteq S(t_2,\kappa_2)$ for $\kappa_1 \geq \kappa_2$. Therefore intervals of the form $S(t,t)$ will be of our interest.

We can take union of $S(t,t)$ with $S(t+1, t+1)$ if $S(t,t) \cap S(t+1, t+1) \neq \varnothing$, which is possible only when the lower bound of $S(t,t)$ is less than the upper bound of $S(t+1, t+1)$. This means
$$
\frac{m t}{t} + \frac{{\sqrt{2m \kappa N_1}}}{\sqrt{t}} \leq \frac{m (t+1)}{t} - \frac{{\sqrt{2m \kappa N_1}}}{\sqrt{t}}.
$$
Thus $t \leq \sqrt{\frac{m}{8N_1}}$, and we get 
$$
\bigcup_{t=1}^{\ceil{\sqrt{\frac{m}{8N_1}}}} S(t,t) \subseteq \left[ m + \sqrt{2mN_1}, \infty \right).
$$
Now $S(t, \kappa) \subseteq U(m,i)$ for $i > 1$, if
$\kappa(i-1)+1 \leq t \leq \kappa i$.
Then for a fixed $\kappa$ we get intervals of the form
$$
S(\kappa(i-1)+1, \kappa), S(\kappa(i-1)+2, \kappa)\cdots S(\kappa(i-1)+\kappa, \kappa). 
$$
Since the lower bound of $S(\kappa(i-1)+2, \kappa)$ is greater than that of $S(\kappa(i-1)+1, \kappa)$ and upper bound of $S(\kappa(i-1)+2, \kappa)$ is less than that of $S(\kappa(i-1)+3, \kappa)$,  so it lies between these two intervals and it it holds for any $S(\kappa(i-1)+j, \kappa)$ for $j= 2,... (\kappa - 1)$. Also, each $S(t, \kappa)$ is contained in some $S(t', \kappa')$ for some $\kappa' < \kappa$. Therefore, like the case $i=1$, here too we find when the boundary intervals for $\kappa$ and $\kappa + 1$ will intersect.

Let us consider
$$
S(t_1, \kappa_1) = S(\kappa(i-1)+1, \kappa) \hspace{5mm} \text{and} \hspace{5mm} S(t_2, \kappa_2) = S((\kappa+1)(i-1)+1, \kappa + 1).
$$
Then again for intersection of them, the upper bound of $S(t_2+1, \kappa_2)$ is $\geq$ the lower bound of $S(t_1, \kappa_1)$ and the bound turns out to be $t_2 < \frac{2m}{16N_1(i-1)^2}$.

Similarly for the other end, taking
$$
S(t_1, \kappa_1) = S(\kappa i, \kappa) \hspace{5mm} \text{and} \hspace{5mm} S(t_2, \kappa_2) = S((\kappa+1) i, \kappa + 1),
$$
we get $t_2 < \frac{m}{16N_1i^2}$.

Taking finite union of $S(t, \kappa)$ considering these two bounds, we get, if
$$
\sqrt{p} \in \left[ \frac{m}{i} + \sqrt{\frac{2mN_1}{i}}, \frac{m}{i-1}- \sqrt{\frac{2mN_1}{(i-1)}}\right],
$$
\end{enumerate}
then there exist totally positive integral elements in $K$ which can not be written as sum of squares.
\end{proof}
We will state some direct consequences of the above theorem.
\begin{cor}\label{cor2}
Let $K = \mathbb{Q}(\sqrt{p},\sqrt{q}) $ where $p , q $ are defined above.
\begin{enumerate}
    \item  If $p,q \geq \left(\frac{m}{2}+4 \right)^2$ then not all elements of $m\mathcal{O}_{F_1}^+$ and $m\mathcal{O}_{F_2}^+$ are represented as the sum of squares in $\mathcal{O}_K$.\label{cor2a}

   \item  If $p \geq \left( m + \sqrt{2mN_1}  \right)^2$ and $q \geq \left( m + \sqrt{2mN_2}  \right)^2$ then not all elements of $m\mathcal{O}_K^+$ are represented as the sum of squares in $\mathcal{O}_K$.\label{cor2b}
\end{enumerate}

\end{cor}

Now the question is, if  $p,q$ are less than the given bound then can we get $F_1 = \mathbb{Q}(\sqrt{p})$ and $F_2 = \mathbb{Q}(\sqrt{q})$ such that every element of $m\mathcal{O}_{F_1}^+$ and $m\mathcal{O}_{F_2}^+$ can be represented as the sum of squares in $\mathcal{O}_K$? Here we have constructed an example for Case - I.

\textit{Example:} For $m = 4$ from  \cite{Ra}, Theorem $7$ we can say for $p \in \{ 2,6,10 \}$ and $q \in \{ 3,7,11 \}$ every element of $4\mathcal{O}_{\mathbb{Q}(\sqrt{p})}^+$ and  $4\mathcal{O}_{\mathbb{Q}(\sqrt{q})}^+$ can be written as sum of squares in  $\mathcal{O}_K = \mathcal{O}_{\mathbb{Q}(\sqrt{p},\sqrt{q})}$.

\section{TABLE I}\label{t1} Let $K = \mathbb{Q}(\sqrt{p},\sqrt{q})$ (where $p, q \geq 2$ as before) and for given positive integer $m$;

 $\alpha_1 \rightarrow$ elements in $\mathcal{O}_K$ which can not be written as sum of squares,

 $\alpha_2 \rightarrow$ elements in $\mathcal{O}_K$ such that $\frac{m}{2}\alpha_3$ can not be written as sum of one square,

$\alpha_3 \rightarrow$ elements in $\mathcal{O}_K$ such that $\frac{m}{2}\alpha_2$ can not be written as sum of two squares. 

The following table exhibits examples of elements in $\mathcal{O}_K^+$ which are of the form $\alpha_1$, $\alpha_2$ and $\alpha_3$.

\begin{center}

    \begin{tblr}{ |c|c|c|c|c|c| } 
\hline
$m$ & $p$ & $q$ & $\alpha_1$ & $\alpha_2$ & $ \alpha_3$\\
\hline
\SetCell{2em}{$2$} & $26$ & $31$ & {$6 +\sqrt{26}$ \\  $6+ \sqrt{31}$} & {$12 + \sqrt{26}+ \sqrt{31}$ \\  $ 23 + 2\sqrt{26}+ 2\sqrt{31}$ \\  $33 + 3\sqrt{26}+ 3\sqrt{31}$} &  $12 + \sqrt{26}+ \sqrt{31} $\\\cline{2-6} 
& $102$ & $203$ & {$11+ \sqrt{102}$ \\ $15+ \sqrt{203}$} & {$ \ceil{k\sqrt{102}} + \ceil{k\sqrt{203}} + k\sqrt{102} + k\sqrt{203}$ \\  where $k \in \{1,...17\}$} & $26 + \sqrt{102} + \sqrt{203}$ \\ 
\hline
\SetCell{2em}{$12$} & $102$ & $103$ & {$11+\sqrt{102}$ \\ $11+\sqrt{103}$} & {$ 22 + \sqrt{102} + \sqrt{103}$ \\ $42 + 2\sqrt{102} + 2\sqrt{103}$} \\\cline{2-6}
& $262$ & $355$ & {$17 + \sqrt{262}$ \\ $19 + \sqrt{355}$} & {$\ceil{k \sqrt{262}} +  \ceil{k \sqrt{355}} + k \sqrt{262} + k \sqrt{355}$ \\ where $k \in \{1,...6\}$} \\
\hline
\end{tblr}
\end{center}
\begin{center}
\begin{tblr}{ |c|c|c|c|c|c| } 
\hline
$m$ & $p$ & $q$ & $\alpha_1$ & $\alpha_2$ & $ \alpha_3$\\
\hline
\SetCell{2em}{$27$} & $326$ & $327$ & {$19 + \sqrt{326}$ \\ $37 + 2\sqrt{326}$ \\ $55 + 3\sqrt{326}$ \\$19 + \sqrt{327}$ \\$37 + 2\sqrt{326}$ \\$55 + 3\sqrt{326}$} & {$38 + \sqrt{326} + \sqrt{327}$ \\  $74 + 2\sqrt{326} + 2\sqrt{327}$ \\$110 + 3\sqrt{326} + 3\sqrt{327}$ \\} \\\cline{2-6}
& $382$ & $511$ & {$20 + \sqrt{382}$ \\ $40 + 2\sqrt{382}$ \\$59 + 3\sqrt{382}$ \\ $23 + \sqrt{511}$\\ $46 + 2\sqrt{511}$\\ $68 + 3\sqrt{511}$\\} & {$43+\sqrt{382}+\sqrt{511}$ \\ $86 + 2\sqrt{382}+ 2\sqrt{511}$ \\ $127 + 3\sqrt{382}+ 3\sqrt{511}$ \\ $170 + 4\sqrt{382} + 4\sqrt{511}$ \\}\\
\hline
\end{tblr}
\end{center}

\section{Case II}\label{c2}
Here the biquadratic field $K = \mathbb{Q}(\sqrt{p},\sqrt{q}) $, where 
$$
p\equiv 2,3 \pmod 4 \hspace{5mm} \text{and} \hspace{5mm} q \equiv 1 \pmod 4.
$$ 
For both of the cases the basis of $\mathcal{O}_K$ is
$$
\mathcal{B}_2 = \{1, \sqrt{p}, \frac{\sqrt{q}+1}{2}, \frac{\sqrt{p}+\sqrt{r}}{2} \}.
$$
The following lemma will be used this in the next proposition. Proof of the first part of it is exactly same as Lemma \ref{lem1} and for the second part see [\cite{Ra}, Lemma 4.2].
\begin{lem}\label{lem12}
Let $a_i$, $b_i$ and $c_i$ be non- negative integers satisfying 
$$
\sum_{i} a_ib_i = \kappa  \hspace{5mm} \sum_{i} a_ic_i = \kappa' 
 \hspace{5mm} \text{and} \hspace{5mm} a_i \equiv c_i \pmod 2
$$ 
for  fixed positive integers $\kappa, \kappa'$. Let $D$ be a real number and $t$ a positive integer such that $t \equiv \kappa' \pmod 2$. Then there exist  intervals
 $$
I_t(\kappa)=\begin{cases} \left[ \frac{\kappa^2}{4t(t+1)},\frac{\kappa^2}{4t(t-1)} \right] & \text{if} \hspace{5mm} t >1 \\
\left[ \frac{\kappa^2}{8},\infty \right) & \text{if} \hspace{5mm} t = 1
\end{cases}
$$
and 
$$
J_t(\kappa)=\begin{cases} \left[ \frac{\kappa'^2}{t(t+2)},\frac{\kappa'^2}{t(t-2)} \right] & \text{if} \hspace{5mm} t >1 \\
\left[ \frac{\kappa'^2}{t(t+2)},\infty \right) & \text{if} \hspace{5mm} t = 1,2
\end{cases}
$$
such that 
$$
\sum_{i} a_i^2 + 4Db_i^2 \geq \frac{\kappa^2}{t}+ 4tD \hspace{5mm} \text{if} \hspace{5mm} D \in I_t(\kappa) 
$$  
and  
$$
\sum_{i} a_i^2 + Db_i^2 \geq \frac{\kappa'^2}{t}+ tD \hspace{5mm} \text{if} \hspace{5mm} D \in J_t(\kappa').
$$
\end{lem}
For integers of this type of quadratic fields, parity plays an important role. Therefore it is necessary to prove a second version of Lemma \ref{lem2}.
\begin{lem}\label{lem22}
Let  $a_i,b_i,c_i$,  $1\leq i \leq n$ be non-negative integers such that
$$
\sum a_ib_i = \kappa_1 \hspace{5mm} \sum a_ic_i = \kappa_2 \hspace{5mm} \sum b_ic_i= \frac{\kappa_3}{g} \hspace{5mm} \text{and} \hspace{5mm} a_i \equiv c_i \pmod 2.
$$
Here $\kappa_1, \kappa_2,\kappa_3$ are any three arbitrary but fixed positive integers, and $t_1, t_2$ be any two positive integers such that $t_2 \equiv \kappa_2 \pmod2$.

If 
$$
p\in \left[ \frac{\kappa_1^2}{4t_1(t_1+1)}, \frac{\kappa_1^2}{4t_1^2}, \right] = H_{t_1}(\kappa_1) \hspace{5mm} \text{and} \hspace{5mm} q\in \left[ \frac{\kappa_2^2}{t_2(t_2+2)}, \frac{\kappa_2^2}{t_2^2}\right] = H'_{t_2}(\kappa_2)  
$$
with $t_1, t_2$ are any two positive integers.
Then,
$$
\sum a_i^2+ 4b_i^2p + c_i^2q \geq \frac{\kappa^2}{t}+ 4t_1p + t_2q,
$$
where $t=$ min$\{ t_1,t_2 \}$ and $\kappa=$ max $\{\kappa_1, \kappa_2 \}$.
\end{lem}

\begin{proof}
Let $l_1= \sum b_i^2$ and $l_2= \sum c_i^2$. Then by Cauchy-Schwarz inequality
$$
\sum a_i^2 \geq \frac{\sum (a_ib_i)^2}{\sum b_i^2}= \frac{\kappa_1^2}{l_1} \hspace{5mm} \text{and} \hspace{5mm} \sum a_i^2 \geq \frac{\sum (a_ic_i)^2}{\sum c_i^2}= \frac{\kappa_2^2}{l_2}.
$$
This implies,
$$
\sum a_i^2+ 4b_i^2p + c_i^2q \geq \frac{\kappa^2}{l}+ 4l_1p + l_2q
$$
Where $l=$ min $\{l_1, l_2\}$ and $\kappa = \text{max} \{\kappa_1, \kappa_2 \}$.

Therefore, we need to show that,
\begin{equation}\label{eqn5}
\frac{\kappa_1^2}{l_1}+ 4l_1p+ l_2q \geq \frac{\kappa_1^2}{t_1}+ 4t_1p + t_2q
\end{equation}
and
\begin{equation}\label{eqn44}
\frac{\kappa_2^2}{l_2}+ 4l_1p+ l_2q \geq \frac{\kappa_2^2}{t_2}+ 4t_1 p+ t_2q.
\end{equation}
Now \eqref{eqn5} implies $(l_1-t_1)(4l_1t_1p - \kappa_1^2) + q(l_2-t_2) \geq 0$ and that holds if 
\begin{equation}\label{eqn55}
4t_1(t_1-1)p \leq \kappa_1^2 \leq 4t_1(t_1+1) \hspace{5mm} \text{and} \hspace{5mm} l_2 \geq t_2.
\end{equation}
Similarly, \eqref{eqn44} implies $(l_2-t_2)(l_2t_2q - \kappa_2^2) + 4p(l_1-t_1) \geq 0$.
Now, $t_2 \equiv \kappa_2 \pmod 2$ and $a_i \equiv c_i \pmod 2$ implies $\kappa_2 \equiv l_2 \pmod 2$ and 
that holds if 
\begin{equation}\label{eqn6}
   t_2(t_2-2)q \leq \kappa_2^2 \leq t_2(t_2+2) \hspace{5mm} \text{and} \hspace{5mm} l_1 \geq t_1.
\end{equation}
Both \ref{eqn55} and \ref{eqn6} are true since $p \in H_{t_1}(\kappa_1) $ and $q \in H'_{t_2}(\kappa_2)$.
\end{proof}
In Lemma  \ref{lem12} the inequalities $\frac{\kappa^2}{l}+ 4lD \geq \frac{\kappa^2}{t}+ tD $ and $\frac{\kappa'^2}{l}+ lD \geq \frac{\kappa'^2}{t}+ tD $ holds if and only if $D \in I_t(\kappa)$ and $D \in J_t(\kappa')$. But this is not the case for Lemma \ref{lem22}. The inequalities $(l_1-t_1)(4l_1t_1p - \kappa_1^2) + q(l_2-t_2) \geq 0$ and $(l_2-t_2)(l_2t_2q - \kappa_2^2) + 4p(l_1-t_1) \geq 0$ can be true by many other ways. Therefore for integers in biquadratic fields we can obtain only necessary condition for not being represented as sum of squares.

Although the next result a lot similar to Proposition \ref{prop1}, but due to the change of basis and the added condition on the parity of $a_i$ and $c_i$ there are subtle differences in the proof.
\begin{prop}\label{prop12}
Let 
$$
I_2= \left[\frac{m\kappa}{4t} + \frac{{\sqrt{m}}}{\sqrt{2t}},\frac{m \kappa}{4(t-1)} - \frac{{\sqrt{m}}}{\sqrt{2(t-1)}} \right]
$$
with $t > 1$
and 
$$
I_3= \left[\frac{m\kappa}{t} + \frac{{2\sqrt{m}}}{\sqrt{t}},\frac{m\kappa}{t-2} - \frac{{2\sqrt{m}}}{\sqrt{t-2}} \right]
$$
for $l>2$ be two intervals.

If $m,\kappa \in \mathbb{Z}_{>0}$, $m \leq$ min $\{ \frac{r}{2 \ceil{\kappa \sqrt{r}}}, p, \frac{q}{2} \} $, $m\kappa \equiv t \pmod 2$ and $\sqrt{p} \in I_{2}$, $\sqrt{q} \in I_{3}$, then not all elements of $m\mathcal{O}_K^{+}$ are represented as sum of integral squares. These intervals are non-empty for 
$$
m\geq \frac{8t(t-1)(2t-1+2\sqrt{t(t-1))}}{k^2}.
$$
For $t=1$ and $t=2$ we define the right bound of the intervals $I_1$ and $I_2$ respectively to be $\infty$ and are therefore always non-empty.
\end{prop}

\begin{proof}
An element $\alpha \in \mathcal{O}_K$ (in the considered field $K$) can be written as 
$$
\alpha = \frac{a}{2} + \frac{b}{2}\sqrt{p} + \frac{c}{2}\sqrt{q} + \frac{d}{2}\sqrt{r}
$$
 with $a,b,c,d \in \mathbb{Z}$, $ 2|b-d $ and $2| a-c$.
 
For arbitrary positive integer $\kappa$ consider $ \alpha = \frac{1}{2}(\ceil{\kappa\sqrt{p}} +\kappa\sqrt{p})$ and suppose $m \alpha = \sum_i \alpha_i^2 $ for some $\alpha_i = \frac{1}{2}(a _i+ b_i\sqrt{p} + c_i\sqrt{q} + d_i\sqrt{r}) $ where $b_i \equiv d_i \pmod 2$ and $a_i \equiv c_i \pmod 2$. Now, comparing trace of $m\alpha$ with that of $\sum_i \alpha_i^2$ we get 
$$
4m \ceil{\kappa\sqrt{p}} = a _i^2 + b_i^2 p  + c_i^2 q + d_i^2 r.
$$
(Since $m < \frac{r}{2 \ceil{\kappa \sqrt{r}}} $ and therefore $m \ceil{\kappa\sqrt{p}} <  \frac{r}{4}$.)

Thus $d_i = 0 $ for all $i$. Then $2|b_i$ and we substitute $b_i$ by $2b_i$ and get $2\alpha_i = a _i+ 2b_i\sqrt{p} + c_i\sqrt{q} $ for all $i$.

Without loss of generality, we can choose $a_i \geq 0$ for all $i$ and $b_i \geq 0$ if $a_i=0$. 
For the sake of contradiction suppose 
Otherwise if there exist some $i$ such that $b_i < 0 , c_i < 0 $ and $a_i > 0$, we have
$$
4m \sigma_4(\alpha) \geq \sigma_4(\alpha_i)^2 = a_i^2 + 4b_i^2 p+ c_i^2q + 4a_i(-b_i) \sqrt{p} + 2a_i(-c_i) \sqrt{q} + 4b_ic_ig \sqrt{r} > 4p \geq 4m
$$
for $p \geq m $, where $g= gcd (p,q)$. Which is impossible because $\sigma_4 (\alpha) < 1$. Therefore, both  $b_i$ and $c_i$ can not be negative for all $i$.

Now, let us suppose that there exist $\alpha_i$ with $b_i < 0$ and $c_i > 0$. Then using $\sigma_2 (\alpha)$ instead of $\sigma_4 (\alpha)$ we get similar type of contradiction. Therefore we can assume $b_i \geq 0$ for all $i$

Comparing the irrational parts in the expression of $m \alpha$, we get $m \kappa = 2\sum_i a_i b_i$ and comparing rational parts from Lemma \ref{lem1} entails 
$$ 
2(m \kappa\sqrt{p} + m) > 2m (\kappa\floor{\sqrt{p}} + 1 ) = \sum_i a_i^2 + 4b_i^2 p + c_i^2 q \geq a_i^2 + 4b_i^2 p \geq \frac{m^2\kappa^2}{4t}+ 4tp
$$
for $p \in I_{t}(\frac{m\kappa}{2})$.
This is impossible if $\sqrt{p} \geq \frac{m \kappa}{4t}+ \frac{\sqrt{m}}{2\sqrt{t}}$ or $\sqrt{p} \leq \frac{m \kappa}{4t} - \frac{\sqrt{m}}{2\sqrt{t}}$.

Combining these restrictions on $p \in I_{t}(\frac{m \kappa}{2})$ and $p \in I_{t-1}(\frac{m \kappa}{2})$ we derive, if 
$$
\sqrt{p} \in \left[ \frac{m \kappa}{4t} + \frac{{\sqrt{m}}}{\sqrt{2t}},\frac{m \kappa}{4(t-1)} - \frac{{\sqrt{m}}}{\sqrt{2(t-1)}} \right] = I_2,
$$
then $m\alpha$ can not be represented as the sum of squares.
This interval is non-empty only when  $|I_2| > 0$. This implies 
\begin{equation}\label{eq1}
    m \geq \frac{8t(t-1)(2t-1+ 2 \sqrt{t(t-1)}}{\kappa^2}.
\end{equation}
For $q$ consider $\alpha = \floor{\frac{\kappa\sqrt{q}-\kappa}{2}}+1+ \frac{\kappa\sqrt{q}-\kappa}{2}$ and $\alpha = \sum_i \alpha_i^2$. Proceeding as in the $p$ case, we can assume that $m \alpha = \sum_i \alpha_i^2 $ for some $\alpha_i = \frac{1}{2}(a _i+ b_i\sqrt{p} + c_i\sqrt{q} + d_i\sqrt{r}) $.

Then as in the previous case,
we get $d_i = 0$ and $c_i \geq 0 $ for all $i$ (since $4m\sigma_3(\alpha)$ and $4m\sigma_4(\alpha)$ are both greater than $q \geq 2m$).

 By comparing irrational part we get $\sum_i a_i c_i = m\kappa$. Also, comparing rational part and using Lemma \ref{lem2}, we get
\begin{eqnarray*}
\frac{m \kappa\sqrt{q}}{2}+m &>& m \left(  \ceil{\frac{\kappa\sqrt{q}-\kappa}{2}} + \frac{\kappa}{2} \right)\\
&=& \sum_i \frac{a_i^2 + 4b_i^2 + c_i^2q}{4} \sum_i \frac{a_i^2 + c_i^2q}{4}\\
&\geq& \frac{m^2\kappa^2}{4t} + \frac{t q}{4}
\end{eqnarray*}
for $q \in J(m\kappa)$ where $t \equiv m \pmod{2}$. The rest of the proof is similar and we get
$$
q \in \left[\frac{m \kappa}{t} + \frac{{2\sqrt{m}}}{\sqrt{t}},\frac{m \kappa}{t-2} - \frac{{2\sqrt{m}}}{\sqrt{t-2}} \right] = I_3
$$
for $t > 2 $.
This is non-empty for 
\begin{equation}\label{eq2}
m \geq  \frac{t(t-2)(2t-2+ 2\sqrt{t(t-2)}}{\kappa^2}.
\end{equation}
And we consider the case $t \in {1,2}$ similarly as in first part.

We get from (\ref{eq1}) and (\ref{eq2}) that all these intervals are nonempty for
$$
m \geq \frac{8t(t-1)(2t-1 + 2 \sqrt{t(t-1)}}{\kappa^2}.
$$
\end{proof}
The next corollary is exactly same as corollary \ref{cor1} and we give the statement as $\alpha_1$ and $\alpha_2$ are not the same. 
\begin{cor}
    Let $\alpha_1 = \frac{1}{2}(\ceil{\kappa\sqrt{p}} +\kappa\sqrt{p}) $, $\alpha_2 = \ceil{\frac{\kappa\sqrt{q} - \kappa}{2}} + \frac{\kappa\sqrt{q} + \kappa}{2} $ and $\alpha_3 =\ceil{\kappa\sqrt{r}} +\kappa\sqrt{r} $ 
    and $m$ is as in Proposition \ref{prop12}. If $m$ is even then
\begin{equation}\label{ka}
\frac{m}{2}(\alpha_i + \alpha_j)
\end{equation}
(where $i,j \in \{1,2,3 \}$ and $i \neq j$) can not be represented as a square in $\mathcal{O}_K$.

If $m =2$ and $\kappa = 1$ then \eqref{ka} can not be written as sum of two integral squares.
\end{cor}
We do not know whether it is possible to represent this type of elements by any number of integer square or not in $K$. Using the method of the above corollary it is difficult to find the same for more than two squares. 

Next proposition is again similar to Proposition \ref{prop12}. Therefore we will only mention some necessary points and skip most of the details.
\begin{prop}\label{prop22}
Let  $ 2\sqrt p + \sqrt{q} + 1 > 2\sqrt{r}$ and $m,t_1,t_2 \in \mathbb{N}$. Let $ \kappa$ be any arbitrary but fixed positive integer such that $2m\kappa \equiv t_2 \pmod 2$ and $m < \frac{\sqrt{r}}{14\kappa}$. Let 
$$
C_1 = \sqrt{q} + 2 \sqrt{r}
$$
and
$$
C_2 = \sqrt{p} - \frac{\sqrt{q}}{2} + \sqrt{r}
$$
be two constants. let
$$
I_3' = \left[ \frac{m \kappa}{4t_1} + \frac{{\sqrt{m\kappa N_1}}}{\sqrt{t_1}},\frac{m \kappa}{4(t_1-1)} - \frac{{\sqrt{m\kappa N_1}}}{\sqrt{t_1-1}} \right] 
$$
for $t_1 >1$, and
$$
I'_4= \left[\frac{m\kappa}{t_2} + 2\frac{{\sqrt{m\kappa N_2}}}{\sqrt{t_2}},\frac{m \kappa}{(t_2-2)} - 2\frac{{\sqrt{m \kappa N_2}}}{\sqrt{t_2-2}} \right]
$$
for $t_2 > 2$ be two intervals. Here $N_1$ and $N_2$ are two constants depending on $C_1$ and $C_2$ respectively.

If $\sqrt{p} \in I'_{3}$ and $\sqrt{q}\in I'_{4}$ then not all elements of $m\mathcal{O}_K^{+}$ are represented as sum of integral squares. These intervals are non-empty for 
$$
    m\geq \text{max} \{ \frac{2N_1t_1(t_1-1)(2t_1-1+2\sqrt{t_1(t_1-1))}}{\kappa}, \frac{8N_2t_2(t_2-2)(t_2-1+\sqrt{t_2(t_2-2))}}{\kappa} \}.
$$
For $t_1 = 1$ and $t_2 = \{1,2\}$ we define the right bound of the intervals to be $\infty$ and are therefore always non-empty.
\end{prop}  
\begin{proof}
    For arbitrary positive integer  $\kappa$ let $\alpha = u + \kappa \sqrt{p} + \frac{\kappa}{2} (1+ \sqrt{q}) + \kappa \sqrt{r}$ where $u = \ceil{\kappa(\sqrt{p} + \frac{1}{2}(\sqrt{q} + 1) + \sqrt{r})}$ and suppose $m \alpha = \sum_i \alpha_i^2 $ for some $\alpha_i = a _i+ \frac{b_i}{2}\sqrt{p} + c_i\sqrt{q} + \frac{d_i}{2}\sqrt{r} $ where $b_i \equiv d_i \pmod 2$. As $m < \frac{\sqrt{r}}{14 \kappa} \implies mu < \frac{r}{4} \implies d_i = 0$ for all $i$. This also implies, $b_i \equiv 0 \pmod 2$. Therefore we substitute $b_i$ by $2b_i$ and get $4\alpha_i = a _i+ 2 b_i\sqrt{p} + c_i\sqrt{q} $ for all $i$.
    
   Then proceeding similarly as in Proposition \ref{prop2}, we can assume $c_i \geq 0$ and $d_i \geq 0$ for all $i$. Now comparing the irrational parts in $m \alpha$, we get,
$$
m \kappa = \sum_i a_i b_i = g\sum_i b_ic_i \hspace{5mm} \text{and} \hspace{5mm } 2m \kappa = \sum_i a_i c_i.
$$
Taking  $N_1 = \frac{1}{2}(C_1 + 3)$ and $N_2 = C_2 + 1$ and comparing the rational parts from Lemma \ref{lem1} we get 
$$ 
4m \kappa(\sqrt{p} + N_1) > mu = \sum_i a_i^2 + b_i^2 p + c_i^2 q \geq \frac{m^2\kappa^2}{t_1}+ 4t_1p + t_2q \geq \frac{m^2\kappa^2}{t_1}+ 4t_1p
$$
for $p \in H'_{t_1}\left( m\kappa \right)$ and 
$$ 
4m \kappa(\sqrt{q} + N_2) > mu \geq \frac{4m^2\kappa^2}{t_2}+ t_2q
$$
for $p \in H'_{t_1}\left( 2m\kappa \right)$.

Combining these restrictions on $p \in H_{t_1}(m \kappa)$ and $p \in H_{t_1-1}(m \kappa)$ we get that, if 
$$
\sqrt{p} \in \left[ \frac{m \kappa}{4t_1} + \frac{{\sqrt{m\kappa N_1}}}{\sqrt{t_1}},\frac{m \kappa}{4(t_1-1)} - \frac{{\sqrt{m\kappa N_1}}}{\sqrt{t_1-1}} \right] 
$$
for $t_1 >1$.
Similarly combining these restrictions for $q \in H'_{t_1}(2m \kappa)$ and $q \in H'_{t_1-2}(2m \kappa)$ we get that, if 
$$
\sqrt{q} \in \left[\frac{m\kappa}{t_2} + 2\frac{{\sqrt{m\kappa N_2}}}{\sqrt{t_2}},\frac{m \kappa}{(t_2-2)} - 2\frac{{\sqrt{m \kappa N_2}}}{\sqrt{t_2-2}} \right]
$$
for $t_2 > 2$,
then $m\alpha$ can not be represented as sum of squares.
Both $I_3'$ and $I_4'$ are non-empty if 
$$
    m\geq \text{max} \{ \frac{4N_it_i(t_i-1)(2t_i-1+2\sqrt{t_i(t_i-1))}}{\kappa} | i= 1,2\}.
$$
Also $I_3'$ and $I_4'$ are well-defined for $t_1, t_2 > 1$. For $t_1 = t_2 = 1$ we get
$$
\sqrt{p} \in \left[\frac{m\kappa}{2} + \sqrt{m\kappa N_1}, \infty \right) \hspace{5mm} \text{for} \hspace{5mm} p \in H_1(m\kappa)
$$
and 
$$
\sqrt{q} \in \left[2m\kappa+ 2\sqrt{m\kappa N_2}, \infty \right) \hspace{5mm} \text{for} \hspace{5mm} p \in H_1'(2m\kappa).
$$
\end{proof}
In Proposition \ref{prop2} and \ref{prop22} we obtain some intervals such that if $\sqrt{p}$ and $\sqrt{q}$ lies in those intervals then the given elements can not be written as sum of integral squares. It is important to observe that it becomes possible to construct such intervals for integers of given form because of some constant relation between $p$ and $q$.

As before, we now try to find the union of the above intervals as in Theorem \ref{thm1}.
\begin{thm}\label{thm2}
Let $m$ be a positive integer and $F_1 = \mathbb{Q}(\sqrt{p})$, $F_2 = \mathbb{Q}(\sqrt{q})$ where $p,q,m$ are defined above in Proposition \ref{prop12} and \ref{prop22}.
\begin{enumerate}
    \item  If 
$$\sqrt{p} \in \left[ \frac{m}{4}+4, \infty \right), \left[ \frac{m}{4i}+ i\sqrt{40}, \frac{m}{4(i-1)} - (i-1)\sqrt{70} \right]
$$
 for integer $i > 1$ and
$$\sqrt{q} \in \left[ \frac{m}{2}+8, \infty \right), \left[ \frac{m}{2i}+ i(2\sqrt{40}), \frac{m}{2(i-1)} - (i-1)(2\sqrt{70}) \right]
$$
for $i > 1$ and $m$ is even, or
$$
\sqrt{q} \in \left[ m+4, \infty \right), \left[ \frac{m}{2i+1}+ (2i+1)(2\sqrt{40}), \frac{m}{2i-1} - (i-1)(2\sqrt{70}) \right]
$$
for $i > 0$ and $m$ is odd,
 then there exist elements in $m\mathcal{O}_{F_1}^{+}$ and $m\mathcal{O}_{F_2}^{+}$ which can not be written as sum of squares in $\mathcal{O}_K$.\label{thm2a}

    \item  If $\sqrt{p} \in \left[ \frac{m}{4} + \sqrt{m N_1}, \infty \right)$, $ \left[ \frac{m}{4i} + \sqrt{\frac{2mN_1}{i}}, \frac{m}{4(i-1)}- \sqrt{\frac{2mN_1}{(i-1)}}\right]$ for integer $i > 1$ and $\sqrt{q} \in \left[ \frac{m}{2} + \sqrt{2m N_2}, \infty \right)$, $ \left[ \frac{m}{2i} + \sqrt{\frac{2mN_2}{i}}, \frac{m}{2(i-1)}- \sqrt{\frac{2mN_2}{(i-1)}}\right]$ for integer $i > 1$. \label{thm2b}

    then there exist elements of $m\mathcal{O}_K^+$ (but not elements of any quadratic subfield of it) which can not be written as sum of squares in $\mathcal{O}_K$.
\end{enumerate} 
\end{thm}
\begin{proof}
As the proof is similar to [\cite{Ra}, Theorem 5], we only derive some upper bounds of those '$t$'s  for which we can take union of the intervals of the form
Let $S(t , \kappa) = \left[ \frac{m \kappa}{4t} + \frac{{\sqrt{m }}}{2\sqrt{t}},\frac{m \kappa}{4(t-1)} - \frac{{\sqrt{m}}}{\sqrt{2(t-1)}} \right]$.
We need to consider only the cases $\frac{t}{\kappa} = 1$ and $\frac{t}{\kappa} > 1$.

For $\frac{t}{\kappa} = 1$ we can take union of $S(t,t)$ with $S(t+1,t+1)$ only when the lower bound of $S(t,t)$ is less than the upper bound of $S(t+1, t+1)$. Which gives $t \leq \frac{m}{32}$ and we get 
$$
\bigcup_{t=1}^{\ceil{\frac{m}{32}}} S(t,t) \subseteq \left[ \frac{m}{4} + 4, \infty \right).
$$
For $\frac{t}{\kappa} > 1$ and for a fixed $\kappa$ we get only $\kappa$ number of intervals of the form 
$$
S(\kappa(i-1)+1, \kappa), S(\kappa(i-1)+2, \kappa),\cdots,  S(\kappa(i-1)+\kappa, \kappa).
$$
Then again for overlapping the boundary interval
$$
S(t_1, \kappa_1) = S(\kappa(i-1)+1, \kappa) \hspace{5mm} \text{with} \hspace{5mm} S(t_2, \kappa_2) = S((\kappa+1)(i-1)+1, \kappa + 1)
$$
we need $t_2 < \frac{m}{2}\frac{1}{4(i-1)^2(3+2\sqrt{2})}$, which gives the upper bound estimate  
$$
\frac{m}{4(i-1)} - (i-1)\sqrt{70}.
$$
Similarly for the other end, taking
$$
S(t_1, \kappa_1) = S(\kappa i, \kappa) \hspace{5mm} \text{and} \hspace{5mm} S(t_2, \kappa_2) = S((\kappa+1) i, \kappa + 1),
$$
we get $t_2 < \frac{m}{2}\frac{1}{4i^2(\frac{5}{2}+\sqrt{6})}$ which gives the lower estimate 
$$
\frac{m}{4i}+ i\sqrt{40}.
$$
For other intervals we can prove the results similarly.
\end{proof}
We note down some of the consequences of the last theorem.
\begin{cor}\label{cor3}
Let $K = \mathbb{Q}(\sqrt{p},\sqrt{q}) $ where $p , q \geq 2$ are defined before.
\begin{enumerate}
    \item  If $p \geq \left(\frac{m}{2}+4 \right)^2$ and $q \geq \left(\frac{m}{2}+8 \right)^2$ and $m$ is even then there exist elements in $m\mathcal{O}_{F_1}^+$ and $m\mathcal{O}_{F_2}^+$ which can not be represented as the sum of squares in $\mathcal{O}_K$. \label{cor3a}

    \item If $p \geq \left(\frac{m}{2}+4 \right)^2$ and $p \geq \left(m + 4 \right)^2$ and $m$ is odd then there exist elements in $m\mathcal{O}_{F_1}^+$ and $m\mathcal{O}_{F_2}^+$ which can not be represented as the sum of squares in $\mathcal{O}_K$.\label{cor3b}

   \item  If $p \geq \left( \frac{m}{4} + \sqrt{mN_1}  \right)^2$ and $q \geq \left( \frac{m}{2} + \sqrt{2mN_2}  \right)^2$ then there exist elements in $m\mathcal{O}_{K}^+$  which can not be represented as the sum of squares in $\mathcal{O}_K$.\label{cor3c}

\end{enumerate}

\end{cor}

Here again for Case II we have constructed examples of $F_1$ and $F_2$ such that every element of $m\mathcal{O}_{F_1}^+$ and $m\mathcal{O}_{F_2}^+$ can be represented as sum of integral squares.

\textit{Example:} For $m = 4$ from  \cite{Ra}, Theorem $7$ we can say for $p \in \{ 2,3,6,7,10,11 \}$ and $q \in \{ 5, 13 \}$ every element of $4\mathcal{O}_{\mathbb{Q}(\sqrt{p})}^+$ and  $4\mathcal{O}_{\mathbb{Q}(\sqrt{q})}^+$ can be written as sum of squares in  $\mathcal{O}_K = \mathcal{O}_{\mathbb{Q}(\sqrt{p},\sqrt{q})}$.
\section{TABLE II}
The following table exhibits some elements in $\mathcal{O}_K^+$ where $K$ is a bi-quadratic field defined in CASE II \eqref{c2} and  $\alpha_1$, $\alpha_2$ and $\alpha_3$ are as in TABLE I \ref{t1}.

\begin{center}

    \begin{tblr}{ |c|c|c|c|c|c| } 
\hline
$m$ & $p$ & $q$ & $\alpha_1$ & $\alpha_2$ & $ \alpha_3$\\
\hline
\SetCell{2em}{$2$} & $23$ & $85$ & {$5 +\sqrt{23}$ \\  $10+ \sqrt{85}$} & {$\ceil{k\sqrt{23}}+\ceil{k\sqrt{85}}+k\sqrt{23}+k\sqrt{85}$ \\ where $k \in \{1,...8\}$ } &  $15 + \sqrt{23}+ \sqrt{85} $\\\cline{2-6} 
& $103$ & $201$ & {$11+ \sqrt{103}$ \\ $15+ \sqrt{201}$} & {$ \ceil{k\sqrt{103}} + \ceil{k\sqrt{201}} + k\sqrt{103} + k\sqrt{201}$ \\  where $k \in \{1,...35\}$} & $26 + \sqrt{103} + \sqrt{201}$ \\ 
\hline
 \SetCell{2em}{$24$} & $106$ & $205$ & {$11 +\sqrt{106}$ \\  $21+ 2\sqrt{106}$ \\ $15 +\sqrt{205}$ \\ $29 +2\sqrt{205}$} & {$ 36 + \sqrt{106} + \sqrt{205}$ \\ $ 50 + 2\sqrt{106} + 2\sqrt{205}$ \\ $ 74 + 3\sqrt{106} + 3\sqrt{205}$  } \\\cline{2-6} 
& $263$ & $453$ & {$17+ \sqrt{263}$ \\ $33 + 2\sqrt{263}$\\$22 + \sqrt{453}$ \\ $43 + 2\sqrt{453}$ \\ } & {$ \ceil{k\sqrt{263}} + \ceil{k\sqrt{453}} + k\sqrt{263} + k\sqrt{453}$ \\  where $k \in \{1,...7\}$}\\ 
\hline  
\end{tblr}
\end{center}
\begin{center}
    \begin{tblr}{ |c|c|c|c|c|c| } 
\hline
$m$ & $p$ & $q$ & $\alpha_1$ & $\alpha_2$ & $ \alpha_3$\\
\hline
\SetCell{2em}{$53$} & $302$ & $1193$ & {$18 +\sqrt{302}$ \\  $35 + \sqrt{1193}$ \\ $35 +2\sqrt{302}$ \\  $70 + 2\sqrt{1193}$ \\$53 +3\sqrt{302}$ \\  $104 + 3\sqrt{1193}$ \\} & {$ 53 + \sqrt{302}+ \sqrt{1193}$ \\ $ 105 + 2\sqrt{302}+ 2\sqrt{1193}$ \\$ 157 + 3\sqrt{302} +3 \sqrt{1193}$ \\$ 209 + 4\sqrt{302} + 4\sqrt{1193}$ \\$ 260 + 5\sqrt{302} +5\sqrt{1193}$ } \\\cline{2-6} 
& $434$ & $1613$ & {$21 + \sqrt{434}$ \\ $41 + \sqrt{1613}$ \\ $42 + 2\sqrt{434}$ \\ $81 + 2\sqrt{1613}$ \\ $63 + 3\sqrt{434}$ \\ $121 + 3\sqrt{1613}$ \\} & {$ \ceil{k\sqrt{434}} + \ceil{k\sqrt{1613}} + k\sqrt{434} + k\sqrt{1613}$ \\  where $k \in \{1,...7\}$} \\ 
\hline    
\end{tblr}
\end{center}
\section{Case III and IV}\label{c3}
In this section we consider real biquadratic fields $K_1$ and $K_2$ both of which are of the form $\mathbb{Q}(\sqrt{p},\sqrt{q}) $ with the following restrictions on $p$ and $q$.
\begin{itemize}
    \item[(i)]  For $K_1$, $ p \equiv 1 \pmod 4 \hspace{5mm} p \equiv 1 \pmod 4 \hspace{5mm} \text{and} \hspace{5mm} gcd(p,q) \equiv 1 \pmod 4  $ with basis 
$$
\mathcal{B}_3 = \{1, \frac{1+ \sqrt{p}}{2}, \frac{\sqrt{q} + 1}{2}, \frac{1 + \sqrt{p} + \sqrt{q} + \sqrt{r}}{4} \}.
$$
\item[(ii)]  For $K_2$, $ p \equiv 1 \pmod 4 \hspace{5mm} p \equiv 1 \pmod 4 \hspace{5mm} \text{and} \hspace{5mm} gcd(p,q) \equiv 3 \pmod 4  $ with basis 
$$
\mathcal{B}_4 = \{1, \frac{1+ \sqrt{p}}{2}, \frac{\sqrt{q} + 1}{2}, \frac{1 - \sqrt{p} + \sqrt{q} + \sqrt{r}}{4} \}.
$$

\end{itemize}
The next result is similar to Proposition \ref{prop1} and Proposition \ref{prop12} but since the basis of $\mathcal{O}_K$ of these fields are different from the first two cases therefore we get some new parity on the coefficients of basis elements which leads to make some relevant changes in the proof. 
\begin{prop}\label{prop61}
Let  $m,t\in \mathbb{N}$ and
$$
I_4= \left[\frac{m\kappa}{t} + 2\frac{{\sqrt{m}}}{\sqrt{t}},\frac{m \kappa}{(t-1)} - 2\frac{{\sqrt{m}}}{\sqrt{t-1}} \right]
$$
for $t > 1$.

If $\kappa \in \mathbb{Z}_{>0}$,  $m \leq$ min $\{ \frac{r}{16 (\kappa \sqrt{r} + 1)}, p, q \} $ and $\sqrt{p}, \sqrt{q} \in I_{4}$ then not all elements of $m\mathcal{O}_{K_{1}^{+}}$ and $m\mathcal{O}_{K_{1}^{+}}$ are represented as sum of integral squares. These intervals are non-empty for 
    $$
    m\geq \frac{t(t-1)(2t-1+2\sqrt{t(t-1))}}{k^2},
    $$
and for $t=1$ we define the right bound of the intervals to be $\infty$ and are therefore always non-empty.
\end{prop}

\begin{proof}
Any element $\alpha \in \mathcal{O}_{K_1}$ can be written as 
$$
\alpha = \frac{1}{4} (a + b\sqrt{p} + c\sqrt{q} +  d\sqrt{r})
$$
 with $a,b,c,d \in \mathbb{Z}$, $ a \equiv b \equiv c \equiv d \pmod 2$ and $a + b + c + d \equiv 0 \pmod 4$.

On the other hand, any element $\alpha \in \mathcal{O}_{K_2}$ can be written as 
$$
\alpha = \frac{1}{4} (a + b\sqrt{p} + c\sqrt{q} +  d\sqrt{r})
$$
 with $a,b,c,d \in \mathbb{Z}$, $ a \equiv b \equiv c \equiv d \pmod 2$ and $a + b + c + d \equiv 0 \pmod 2$.

For arbitrary positive integer $\kappa$ consider $ \alpha = \ceil{\frac{\kappa\sqrt{p} - \kappa}{2}} +\frac{\kappa + \kappa\sqrt{p}}{2} \in \mathcal{O}^+ _{\mathbb{Q}(\sqrt{p})}$ from both $K_1$ and $K_2$. Suppose $m \alpha = \sum_i \alpha_i^2 $ for some $\alpha_i = \frac{1}{4}(a _i+ b_i\sqrt{p} + c_i\sqrt{q} + d_i\sqrt{r})$ where $a_i, b_i, c_i, d_i$ are one of the above type.

Since $m < \frac{r}{16 (\kappa \sqrt{r} + 1)} $ implies $d_i = 0$, from the given conditions we get for both of the fields that  $a_i, b_i, c_i$ are even for all $i$. Therefore we substitute $a_i, b_i, c_i$ by $2a_i, 2b_i$ and $2 c_i$ respectively.
Remaining of the proof is similar to Proposition \ref{prop1}. 
\end{proof}
\begin{cor}
    Let $\alpha_1 = \ceil{\frac{\kappa\sqrt{p} - \kappa}{2}} +\frac{\kappa + \kappa\sqrt{p}}{2} $, $\alpha_2 = \ceil{\frac{\kappa\sqrt{q} - \kappa}{2}} +\frac{\kappa + \kappa\sqrt{q}}{2} $ and $\alpha_3 = \ceil{\frac{\kappa\sqrt{r} - \kappa}{2}} +\frac{\kappa + \kappa\sqrt{r}}{2} $ where 
    $ m$ is defined in Proposition \ref{prop61}. If $m$ is even then 
\begin{equation}\label{ka1}
\frac{m}{2}(\alpha_i + \alpha_j)
\end{equation}
(where $i,j \in \{1,2,3 \}$ and $i \neq j$)
can not be represented as a square in $\mathcal{O}_K$.

If $m =2$ and $\kappa = 1$ then \eqref{ka1} can not be written as sum of two integral squares.
\end{cor}
In the next Lemma Minkowski's inequality will be used in a different way so that we can get same interval $E_{t_1+t_2}(\kappa_1 + \kappa_2)$ for both $p$ and $q$.

\begin{lem}\label{lem31}
Let $p < q < r$ be as before.
Also let $a_i,b_i,c_i$,  $1\leq i \leq n$ be non-negative integers satisfying
$$
\sum a_ib_i = \kappa_1 \hspace{5mm} \sum a_ic_i = k_2 \hspace{5mm} \sum b_ic_i= \frac{\kappa_3}{g}
$$

where $\kappa_1, \kappa_2,\kappa_3$ are any three fixed positive integers and $g = gcd(p,q)$.

If 
$$
p , q\in \left[ \frac{(\kappa_1 + \kappa_2)^2}{2(t_1+t_2)(t_1+t_2+1)}, \frac{(\kappa_1 + \kappa_2)^2}{2(t_1+t_2)^2}\right] = E_{t_1+t_2}(\kappa_1 + \kappa_2)   
$$
where $t_1, t_2$ are any two positive integers, then
$$
\sum a_i^2+ b_i^2p + c_i^2q \geq \frac{(\kappa_1 + \kappa_2)^2}{2(t_1 + t_2)}+ t_1p + t_2q.
$$
\end{lem}

\begin{proof}
Let $l_1= \sum b_i^2$ and $l_2= \sum c_i^2$. Then by Cauchy-Schwarz inequality
\begin{equation} \label{eqn61}
    \sum a_i^2 \sum b_i^2 \geq \sum (a_i b_i)^2 = \kappa_1^2 \hspace{5mm} \implies \sum a_i^2 l_1 \geq \kappa_1^2 
\end{equation}|
and
\begin{equation} \label{eqn62}
\sum a_i^2 \sum c_i^2 \geq \sum (a_i c_i)^2 = \kappa_2^2 \hspace{5mm} \implies \sum a_i^2 l_2 \geq \kappa_2^2.
\end{equation}
From \eqref{eqn61} and \eqref{eqn62} we get 
$$
\sum a_i^2 \geq \frac{\kappa_1^2 + \kappa_2^2}{l_1 + l_2} \geq \frac{(\kappa_1 + \kappa_2)^2}{2(l_1 + l_2)}.
$$
This implies,
$$
\sum a_i^2+ b_i^2p + c_i^2q \geq \frac{(\kappa_1 + \kappa_2)^2}{2(l_1 + l_2)} + l_1p + l_2q.
$$
Therefore, we have to show that,
\begin{equation}\label{eqn 63}
\frac{(\kappa_1 + \kappa_2)^2}{2(l_1 + l_2)} + l_1p + l_2q \geq \frac{(\kappa_1 + \kappa_2)^2}{2(t_1 + t_2)} + t_1p + t_2q.
\end{equation}
Now \eqref{eqn 63} implies 
$$
(l_1-t_1)\left((l_1 + l_2)(t_1 + t_2)p - \frac{(\kappa_1 + \kappa_2)^2}{2}\right) + (l_2-t_2)\left((l_1 + l_2)(t_1 + t_2)q - \frac{(\kappa_1 + \kappa_2)^2}{2}\right) \geq 0,
$$
which is true if 
\begin{equation}\label{eqn64}
(l_1-t_1)\left((l_1 + l_2)(t_1 + t_2)p - \frac{(\kappa_1 + \kappa_2)^2}{2}\right) \geq 0
\end{equation}
and 
\begin{equation}\label{eqn65}
    (l_2-t_2)\left((l_1 + l_2)(t_1 + t_2)q - \frac{(\kappa_1 + \kappa_2)^2}{2}\right) \geq 0.
\end{equation}
Further \eqref{eqn64} implies 
\begin{equation}\label{eqn 66}
  (l_1-t_1) \geq 0 \hspace{5mm} \text{and} \hspace{5mm} \left((l_1 + l_2)(t_1 + t_2)p - \frac{(\kappa_1 + \kappa_2)^2}{2}\right) \geq 0
\end{equation}
Or,
\begin{equation}\label{eqn 67}
  (l_1-t_1) \leq 0 \hspace{5mm} \text{and} \hspace{5mm} \left((l_1 + l_2)(t_1 + t_2)p - \frac{(\kappa_1 + \kappa_2)^2}{2}\right) \leq 0.
\end{equation}
 Similarly from \eqref{eqn65} we get 
 \begin{equation}\label{eqn68}
    (l_2-t_2) \geq 0 \hspace{5mm} \text{and} \hspace{5mm} \left((l_1 + l_2)(t_1 + t_2)q - \frac{(\kappa_1 + \kappa_2)^2}{2}\right) \geq 0 
 \end{equation}
Or,
\begin{equation}\label{eqn69}
    (l_2-t_2) \leq 0 \hspace{5mm} \text{and} \hspace{5mm} \left((l_1 + l_2)(t_1 + t_2)q - \frac{(\kappa_1 + \kappa_2)^2}{2}\right) \leq 0.
 \end{equation}
One of \eqref{eqn68} or \eqref{eqn69} must be true since 
$$
p \in \left[ \frac{(\kappa_1 + \kappa_2)^2}{2(t_1+t_2)(t_1+t_2+1)}, \frac{(\kappa_1 + \kappa_2)^2}{2(t_1+t_2)^2}\right] = E_{t_1+t_2}(\kappa_1 + \kappa_2).  
$$
Similarly we get the upper and lower bounds for $q$. 
\end{proof}
Next proposition is similar to proposition \ref{prop2} and proposition \ref{prop22} but since the interval $E_{t_1+t_2}(\kappa_1 + \kappa_2)$ is not similar to $L_t(\kappa),L'_t(\kappa'),H_t(\kappa),H'_t(\kappa'),$ therefore we will make some necessary changes in the proof here.

\begin{prop}\label{prop62}
    Let $ \sqrt p + \sqrt{q} + 1 > \sqrt{r}$ and $m,t_1,t_2 \in \mathbb{N}$ with $t_1 \leq t_2$. Let $ \kappa$ be any arbitrary but fixed even positive integers for $K_1$ and any arbitrary but fixed positive integer (need not be even) for $K_2$ and $m < \frac{\sqrt{r}}{24\kappa}$. Let us also assume 
$$
C_1 = \sqrt{q} + \sqrt{r}
$$
And
$$
C_2 = \sqrt{p} + \sqrt{r}
$$
are two constants.
Let 
$$
I_5' = \left[ \frac{m \kappa}{2(t_1 + t_2)} + 2\frac{{\sqrt{m\kappa N_1}}}{\sqrt{t_1 + t_2}},\frac{m \kappa}{2(t_1 + t_2)-1} - 2\frac{{\sqrt{m\kappa N_1}}}{\sqrt{(t_1 + t_2)-1}} \right] 
$$
$$
I'_6= \left[\frac{m\kappa}{t_1 + t_2} + \frac{{\sqrt{8m\kappa N_2}}}{\sqrt{t_1 + t_2}},\frac{m \kappa}{(t_1 + t_2)-1} - \frac{{\sqrt{8m \kappa N_2}}}{\sqrt{(t_1 + t_2)-1}} \right]
$$
be two intervals, where $N_1$ and $N_2$ are two constants depending on $C_1$ and $C_2$ respectively.

If $\sqrt{p} \in I'_{5}$ and $\sqrt{q} \in I'_{6}$ then not all elements of $m\mathcal{O}_K^{+}$ are represented as sum of integral squares. The interval $I'_5$ is non-empty for 
$$
m\geq   \frac{16N_1t(t-1)(2t-1+2\sqrt{t(t-1))}}{\kappa}.
$$
Also $I'_6$ is non-empty for 
$$
m\geq   \frac{8N_2t(t-1)(2t-1+2\sqrt{t(t-1))}}{\kappa}.
$$
Here $t = t_1 + t_2$
\end{prop}
\begin{proof}
For both $K_1$ and $K_2$ we consider elements of the same form and proceed exactly in the same way. As the technique of the proof is similar to Proposition \ref{prop2} and \ref{prop22}, we only mention some essential points.

For arbitrary positive integer  $\kappa$ (even for $K_1$) let $\alpha = u + \kappa \frac{\sqrt{p} + 1}{4} + \kappa \frac{\sqrt{q} + 1}{4}  + \kappa \frac{\sqrt{r} + 1}{4}$ where $u = \ceil{\frac{\kappa(\sqrt{p} + \sqrt{q} + \sqrt{r}) - 3 \kappa}{4}}$. Then proceeding similarly as in Proposition \ref{prop2} and Proposition \ref{prop22} we get $b_i \geq 0$ and $c_i \geq 0$ for all $i$. 

If we take $N_1 = \frac{C_1 + 4}{4}$ and $N_2 = \frac{C_2 + 4}{4}$ then after comparing the rational and irrational parts in the expression of $m \alpha$ we get 

$$ 
4m \kappa(\frac{\sqrt{p}}{4} + N_1) > mu = \sum_i a_i^2 + b_i^2 p + c_i^2 q \geq \frac{m^2\kappa^2}{2(t_1 + t_2)}+ t_1p + t_2q \geq \frac{m^2\kappa^2}{t_1}+ t_1p
$$
for $p \in E_{t_1 + t_2}\left( m\kappa \right)$, and
$$ 
4m \kappa(\frac{\sqrt{q}}{4} + N_1) > mu = \sum_i a_i^2 + b_i^2 p + c_i^2 q \geq \frac{m^2\kappa^2}{2(t_1 + t_2)}+ t_1p + t_2q.
$$
Which implies
$$
8m \kappa(\frac{\sqrt{q}}{4} + N_1) \geq \frac{m^2\kappa^2}{(t_1 + t_2)} + (t_1 + t_2) q \hspace{5mm} \text{(Using $t_1 \leq t_2$)}
$$
For $q \in E_{t_1 + t_2}\left( m\kappa \right)$.
Combining these restrictions for $p,q  \in E_{t_1 + t_2}(m \kappa)$ and $p,q \in E_{(t_1 + t_2)-1}(m \kappa)$, we get that, if 
$$
\sqrt{p} \in \left[ \frac{m \kappa}{2(t_1 + t_2)} + 2\frac{{\sqrt{m\kappa N_1}}}{\sqrt{t_1 + t_2}},\frac{m \kappa}{2(t_1 + t_2)-1} - 2\frac{{\sqrt{m\kappa N_1}}}{\sqrt{(t_1 + t_2)-1}} \right] = I'_5,
$$
and 
$$
 \sqrt{q} \in \left[\frac{m\kappa}{t_1 + t_2} + \frac{{\sqrt{8m\kappa N_2}}}{\sqrt{t_1 + t_2}},\frac{m \kappa}{(t_1 + t_2)-1} - \frac{{\sqrt{8m \kappa N_2}}}{\sqrt{(t_1 + t_2)-1}} \right] = I'_6,
$$
then $m\alpha$ can not be represented as the sum of squares.
 
\end{proof}
Since $t_1$ and $t_2$ are both positive integers, it is important to note that for any value of $t_1$ and $t_2$ the upper bound of the intervals $I'_5$ and $I'_6$ can not be infinity.

The next theorem is similar to Theorem \ref{thm1} and Theorem \ref{thm2}. We will only state it.
\begin{thm}\label{thm3}
Let $m$ be a positive integer and $F_1 = \mathbb{Q}(\sqrt{p})$, $F_2 = \mathbb{Q}(\sqrt{q})$ where $p,q,m$ are defined above in Propositions \ref{prop61} and \ref{prop62}.
\begin{enumerate}
\item  If 
$$
\sqrt{p}, \sqrt{q} \in \left[ m+8, \infty \right), \left[ \frac{m}{i}+ 2i\sqrt{40}, \frac{m}{(i-1)} - 2(i-1)\sqrt{70} \right]
$$

for integer $i > 1$ then there exist elements in $m\mathcal{O}_{F_1}^{+}$ and $m\mathcal{O}_{F_2}^{+}$ which can not be written as sum of squares in $\mathcal{O}_K$. \label{thm3a}

    \item  If 
$$
\sqrt{p} \in \left[ \frac{m}{2} + 2\sqrt{m N_1}, \infty \right),  \left[ \frac{m}{2i} + 2\sqrt{\frac{mN_1}{i}}, \frac{m}{2(i-1)}- 2\sqrt{\frac{mN_1}{(i-1)}}\right]
$$
 for integer $i > 1$ and 
$$
\sqrt{q} \in \left[ m + \sqrt{8m N_2}, \infty \right),  \left[ \frac{m}{i} + \sqrt{\frac{8mN_2}{i}}, \frac{m}{i-1}- \sqrt{\frac{8mN_2}{(i-1)}}\right]
$$
for integer $i > 1$,
then there exist elements of $m\mathcal{O}_K^+$ (but not elements of any quadratic subfield of it) which can not be written as sum of squares in $\mathcal{O}_K$.\label{thm3b}
\end{enumerate} 
\end{thm}
We are concluding this section with the following corollary similar to Corollary \ref{cor2} and Corollary \ref{cor3}
\begin{cor}\label{cor4}
Let $K = \mathbb{Q}(\sqrt{p},\sqrt{q}) $ where $p , q $ are defined above.
\begin{enumerate}
    \item  If $p,q \geq \left(m + 8 \right)^2$ then not all elements of $m\mathcal{O}_{F_1}^+$ and $m\mathcal{O}_{F_2}^+$ are represented as the sum of squares in $\mathcal{O}_K$.\label{cor4a}
   \item  If $p \geq \left( \frac{m}{2} + 2\sqrt{mN_1}  \right)^2$ and $q \geq \left( m + \sqrt{8mN_2}  \right)^2$ then not all elements of $m\mathcal{O}_K^+$ are represented as the sum of squares in $\mathcal{O}_K$.\label{cor4b}
\end{enumerate}
\end{cor}

For CASE III, examples of $F_1$ and $F_2$ in which every totally positive integer can be represented as sum of integer squares are given here.
 
\textit{Example:} For $m = 4$ from  \cite{Ra}, Theorem $7$ we can say for $p \in \{ 5, 13 \}$ and $q \in \{ 5, 13 \}$ every element of $4\mathcal{O}_{\mathbb{Q}(\sqrt{p})}^+$ and  $4\mathcal{O}_{\mathbb{Q}(\sqrt{q})}^+$ can be written as sum of squares in  $\mathcal{O}_K = \mathcal{O}_{\mathbb{Q}(\sqrt{p},\sqrt{q})}$.
 
\section{TABLE III}
The following table provides examples of elements in $\mathcal{O}_K^+$ where $K$ is a bi-quadratic field defined in CASE III and IV \ref{c3} and  $\alpha_1$, $\alpha_2$ and $\alpha_3$ are defined as in TABLE I \ref{t1}.
\begin{center}
    \begin{tblr}{ |c|c|c|c|c|c| } 
\hline
$m$ & $p$ & $q$ & $\alpha_1$ & $\alpha_2$ & $ \alpha_3$\\
\hline
 \SetCell{2em}{$2$} & $101$ & $105$ & {$11 +\sqrt{101}$  \\ $11 +\sqrt{105}$ } & {$ 22 + \sqrt{101} + \sqrt{105}$ \\ $ 42 + 2\sqrt{101} + 2\sqrt{105}$ \\ $ 62 + 3\sqrt{101} + 3\sqrt{105}$  } & $ 22 + \sqrt{101} + \sqrt{105}$ \\\cline{2-6} 
& $113$ & $251$ & {$11+ \sqrt{113}$ \\ $16 + \sqrt{251}$} & {$ \ceil{k\sqrt{113}} + \ceil{k\sqrt{251}} + k\sqrt{113} + k\sqrt{251}$ \\  where $k \in \{1,...5\}$} & $27+\sqrt{113}+\sqrt{251}$\\ 
\hline
\SetCell{2em}{$5$} & $173$ & $177$ & {$14 +\sqrt{173}$ \\  $14 + \sqrt{177}$ \\ $27 +2\sqrt{173}$ \\  $27 + 2\sqrt{177}$ } & {$ 28 + \sqrt{173}+ \sqrt{177}$ \\ $ 54 + 2\sqrt{173}+ 2\sqrt{177}$} \\\cline{2-6} 
& $209$ & $321$ & {$15 + \sqrt{209}$ \\ $18 + \sqrt{321}$ \\ $29 + 2\sqrt{209}$ \\ $36 + 2\sqrt{321}$ } & {$ 33+ \sqrt{209} + \sqrt{321}$ \\$ 65 + 2\sqrt{209} + 2\sqrt{321}$ \\$ 98 + 3\sqrt{209} + 3\sqrt{321}$ \\} \\ 
\hline  
\SetCell{2em}{$7$} & $229$ & $233$ & {$16 +\sqrt{229}$ \\  $16 + \sqrt{233}$ \\ $31 +2\sqrt{229}$ \\  $31 + 2\sqrt{233}$ } & {$ 32 + \sqrt{229}+ \sqrt{233}$ \\ $ 62 + 2\sqrt{229}+ 2\sqrt{233}$} \\\cline{2-6} 
& $249$ & $465$ & {$16 + \sqrt{249}$ \\ $22 + \sqrt{465}$ \\ $32 + 2\sqrt{249}$ \\ $44 + 2\sqrt{465}$ \\ $113 + 3\sqrt{249}$ \\ $65 + 3\sqrt{465}$ } & {$ 38+ \sqrt{249} + \sqrt{465}$ \\$ 76 + 2\sqrt{249} + 2\sqrt{465}$ \\$ 113 + 3\sqrt{249} + 3\sqrt{465}$ \\} \\ 
\hline    
\end{tblr}
\end{center}

\section*{Acknowledgments}

The author is thankful to her adviser Prof. Kalyan Chakraborty  and Dr. Azizul Hoque for introducing her into this beautiful area of research. She is indebted to Prof Vita Kala for many a discussions that she had with him and also for his valuable suggestions which has helped improving the results in the manuscript immensely.


\begin{thebibliography}{25}
\bibitem{LA}  J. L. Lagrange, {\it D\'{e}monstration d'un th\'{e}or\'{e}me d'arithm\'{e}tique}, in: Nouveaux M\'{e}m. Acad. Roy. Sci. Belles-Lettres, Berlin, 1770, reprinted in: \OE uvres {\bf 3} (1869), 189--201.

\bibitem{OM}  O. T. O'mear, Introduction to quadratic forms, Grundlehren Math., vol. 117, Springer-Verlag, 1971.

\bibitem{SI21} C. L. Siegel, {\it Darstellung total positive Zahlen durch Quadrate}, Math. Z. {\bf 11} (1921), 246--275.

\bibitem{SI45} C. L.Siegel, {\it Sums of $m$-th powers of algebraic integers}, Ann. of Math. {\bf 46} (1945, 313--339.

\bibitem{Wi} K. S. Williams, {\it Integers of biquadratic fields}, Canad. Math. Bull., Vol. 13 (1970), 519–526.

\bibitem{KY} V. Kala and P. Yatsyna, {\it Sums of squares in S-integers}, New York J. Math {\bf 26} (2020), 1145-1154.

\bibitem{KY2} V. Kala and P. Yatsyna, { \it Lifting problems for universal quadratic forms}, Adv. Math. {\bf 377} (2021), 24pp

\bibitem{D} L. E. Dickson, {\it Quaternary Quadratic Forms Representing All Integers}, Amer. J. Math. {\bf 49} (1927) , 39-56.

\bibitem{R} S. Ramanujan, {\it On the Expression of a Number in the Form $ax^2+by^2+cz^2+du^2$}, Proc. Cambridge Phil. Soc. {\bf 19} (1917), 11-21.

\bibitem{Ja} F. Javris, {\it Algebraic Number Theory}, Springer, 2007.

\bibitem{Ra} M. Raska, {\it Representing Multiple of $M$ In Real Quadratic Fields As Sums Of Squares}, arXiv:2105.11423v1


\bibitem{Co} J. H. Conway, {\it Universal quadratic forms and the fifteen theorem}, Contemp. Math. {\bf 272} (1999), 23-26


\bibitem{BH} M. Bhargava, J.Hanke, {\it Universal quadratic forms and the $290$ theorem}, preprint

\bibitem{Ro} J.Rouse, {\it Quadratic forms representing all odd positive integers}, Amer. J. Math. {\bf 136} (2014), 1693-1745


\end{thebibliography}
\end{document}